\documentclass[11pt]{pjm} 
  
\usepackage{graphics}

\title[Toda latttice]{ Toda lattice and toric varieties for real split
semisimple Lie algebras } 
  
\author[Casian]{Luis G. Casian} 
\address{ Department of Mathematics \\
The Ohio State University \\
Columbus, OH 43210 .} 
\email{casian@math.ohio-state.edu} 

\author[Kodama]{Yuji Kodama} 
\address{ Department of Mathematics \\
The Ohio State University \\
Columbus, OH 43210 .} 
\email{kodama@math.ohio-state.edu} 

\thanks{ This research was supported by NSF grant DMS-0071523}

\newtheorem{thm}{Theorem}[section] 
\newtheorem{prop}[thm]{Proposition} 
\newtheorem{lem}[thm]{Lemma} 
\newtheorem{cor}[thm]{Corollary} 
  
\theoremstyle{definition} 
  
\newtheorem{defn}[thm]{Definition} 
\newtheorem{notation}[thm]{Notation} 
\newtheorem{example}[thm]{Example}

\theoremstyle{remark} 
  
\newtheorem{rem}[thm]{Remark}

\def\proof{\par{\it Proof}. \ignorespaces}
\def\endproof{{\ \vbox{\hrule\hbox{%
\vrule height1.3ex\hskip0.8ex\vrule}\hrule }}\par}
\newenvironment{Proof}{\proof}{\endproof}

\date{Received date / Revised version date}

\begin{document} 
  
\begin{abstract} 
  
The paper concerns the topology of an isospectral
real smooth manifold for certain Jacobi element associated with real
split semisimple Lie algebra.
The manifold is
identified as a compact, connected completion
of the disconnected Cartan subgroup of the corresponding
Lie group $\tilde G$ which is a disjoint
union of the split Cartan subgroups associated to
semisimple portions of Levi factors of all standard
parabolic subgroups of $\tilde G$.
The manifold is also related to the
compactified level sets of a generalized Toda lattice
equation defined on the semisimple Lie algebra, which is
diffeomorphic to a toric variety in the flag manifold ${\tilde G}/B$
with Borel subgroup $B$ of $\tilde G$.
 We then give a cellular
decomposition and the associated chain complex
of the manifold by introducing colored-signed
Dynkin diagrams which parametrize the cells in the decomposition.
  
\end{abstract} 
  
\maketitle 
  
\section{Introduction}\label{S1}

In this paper, we study the topological structure of certain manifolds that  are interesting in two different ways. First they are isospectral  manifolds for a signed Toda lattice flow \cite{KO-YE1}; an integrable system that arises in several physical  contexts and has been studied extensively. Secondly  they are shown in \S \ref{S8} to be the closures of generic orbits of a split Cartan subgroup on a real  flag manifold.  These are certain smooth toric varieties that glue together the disconnected pieces of a Cartan subgroup of a semisimple Lie group of the form $\tilde G$.  In this paper we start from the Toda lattice aspect of this object and end the paper inside a real flag manifold. We thus start by motivating and describing our main constructions from the point of view of the Toda lattice; then we trace a path that starts with this Toda lattice and that naturally leads to the (disconnected) Cartan subgroup of a split semisimple Lie group and a toric orbit. Another path that also leads to a Cartan subgroup starts with Kostant's paper \cite{KOS}; this approach is described in \S \ref{S8}.

From the Toda lattice end of this story, these manifolds are related to the compactified level set
of a generalized (nonperiodic) Toda lattice equation defined on the semisimple
Lie algebra (see for example \cite{KOS}) and, although they share some features
with the Tomei manifolds in \cite {TO}, they  are different from those (e.g., non-orientable).  As  background information, we
start with a definition of the generalized Toda lattice equation
which led us to our present study of the manifolds.

Let $\mathfrak g$ denote a real
split semisimple Lie algebra of rank $l$. We fix a
split Cartan subalgebra $\mathfrak h$ with root system $\Delta
=\Delta({\mathfrak g},{\mathfrak h}) $, real root vectors
$e_{\alpha_i}$ associated with simple roots $ \left\{ \alpha_i : i=1,..,l \right\} =\Pi$.
We also denote $\left\{ h_{\alpha_i},e_{\pm\alpha_i} \right\} $ the Cartan-Chevalley basis of
$\mathfrak g$ which satisfies the relations,
\begin{equation}
 {[h_{\alpha_i} , h_{\alpha_j}] = 0,} \quad
 { [h_{\alpha_i}, e_{\pm \alpha_j}] = \pm C_{j,i}e_{\pm \alpha_j}~,} \quad
{[e_{\alpha_i} , e_{-\alpha_j}] = \delta_{i,j}h_{\alpha_j}.}
\nonumber
\end{equation}
where the $l\times l$ matrix
$(C_{i,j})$ is the Cartan matrix
corresponding to $\mathfrak g$, and $C_{i,j}=\alpha_i(h_{\alpha_j})$.

 Then the generalized Toda lattice equation related to
real split semisimple Lie algebra is defined by
the following system of 2nd order differential equations for
the real variables $f_i(t)$ for $i=1,\cdots,l$,
\begin{equation}
 {d^2 f_i \over dt^2}= \epsilon_i \exp \left(
-\sum_{j=1}^l C_{i,j}f_j\right)  
\label{1.1.1}
\end{equation}
where $\epsilon_i\in \{\pm 1\}$ which correspond to the signs
in the indefinite Toda lattices introduced in \cite{faybusovich:90, KO-YE1}.
The main feature of the indefinite Toda equation having
at least one of $\epsilon_i$ being $-1$ is that the
solution blows up to infinity in finite time \cite{KO-YE1, gekhtman:97}.
Having introduced the signs, the group corresponding to the Toda lattice is
a real split Lie group $\tilde G$ with Lie algebra
$\mathfrak g$ which is defined in  \S \ref{S2}.
For example, in the case of ${\mathfrak g}={\mathfrak{sl}}(n,\mathbb R)$,
if $n$ is odd, $\tilde
G=SL(n,\mathbb R)$, and if $n$ is even, $\tilde G=Ad( SL(n,\mathbb R)^{\pm})$.

The original Toda lattice equation in \cite{T} describing a system of
$l$ particles on a line
interacting pairwise with exponential forces corresponds to the case with
${\mathfrak g}={\mathfrak{sl}}(l+1,\mathbb R)$ and $\epsilon_i=1$ for all $i$,
and it is given by
\begin{equation}
{d^2 q_i \over dt^2}=  \exp (q_{i-1}-q_{i})-\exp (q_i-q_{i+1}), \quad
i=1,\cdots,l,
\nonumber
\end{equation}
where the physical variable $q_i$, the position
of the $i$-th particle, is given by
\begin{equation}
 \displaystyle{ q_i=f_{i}-f_{i+1}, } \quad
 \displaystyle{ i=1,\cdots,l,}
\nonumber
\end{equation}
with $f_{l+1}=0$ and $f_{0}=f_{l+2}=-\infty$
indicating $q_0=-\infty$ and $q_{l+1}=\infty$.

The system (\ref{1.1.1}) can be written in the so-called Lax equation
which describes an iso-spectral deformation of a Jacobi element of
the algebra $\mathfrak g$.  This is formulated by defining the set of real functions $\{(a_i(t),b_i(t)):i=1,\cdots,l\}$ with
\begin{equation}
\begin{array} {lll}
 &a_i(t):=\displaystyle{{d f_i(t) \over dt},} 
 &b_i(t):=\displaystyle{\epsilon_i \exp \left(-\sum_{j=1}^l C_{i,j}f_j(t)\right) }\\
\label{1.1.2}
\end{array}
\end{equation}
from which the system (\ref{1.1.1}) reads
\begin{equation}
\displaystyle{d a_i \over dt}=b_i, \quad \quad
\displaystyle{{d b_i \over dt}=  - b_i 
\left(\sum_{j=1}^l C_{i,j} a_j\right)} .
\label{1.1.3b}
\end{equation}
This is then equivalent to the Lax equation defined on $\mathfrak g$
(see \cite{FL} for a nice review of the Toda equation).
\begin{equation}
{dX(t) \over dt}=[P(t), X(t)] ,
\label{1.1.4}
\end{equation}
where the Lax pair $(X(t),P(t))$ in ${\mathfrak g}$ is defined by
\begin{equation}
\left\{
\begin{array}{ll}
 X(t) &=\displaystyle{\sum_{i=1}^l a_i(t)h_{\alpha_i}+\sum_{i=1}^l \left(
b_i(t)e_{-\alpha_i}+e_{\alpha_i}\right)} \\
P(t) & =\displaystyle{-\sum_{i=1}^l b_i(t)e_{-\alpha_i}}.
\end{array}
\right.
\label{1.1.5b}
\end{equation}

Although a case with some $b_i=0$ is not defined in (\ref{1.1.2})  the corresponding
system (\ref{1.1.3b}) is well defined and is reduced to several noninteracting
subsystems separated by $b_i=0$. The constant solution $a_i$ for
$b_i=0$ corresponds to an eigenvalue of the Jacobi (tridiagonal) matrix
for $X(t)$ in the adjoint representaion of ${\mathfrak g}$.
We denote by $S(F)$  the $ad$ diagonalizable elements in
${\mathfrak g } {\otimes }_{\mathbb R} F$
 with eigenvalues in $F$,
where $F$ is $\mathbb R$ or $\mathbb C$.

\begin{figure}
\includegraphics{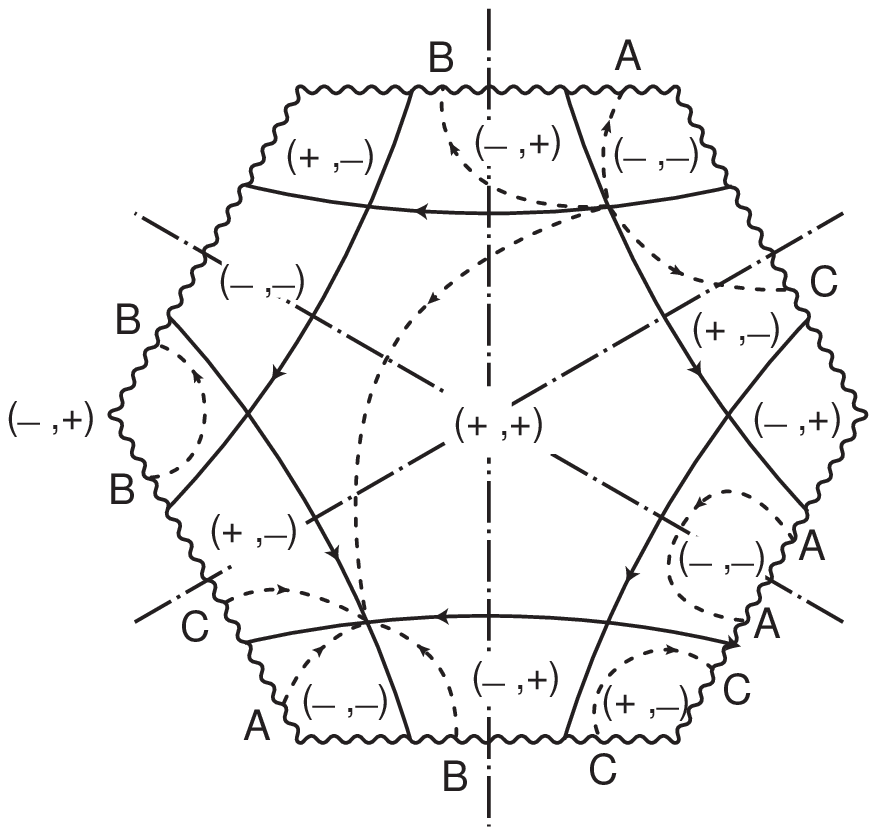}
\caption{The manifold $Z(\gamma)_{\mathbb R}$ and the Toda orbits for ${\mathfrak {sl}}(3, \mathbb R)$.}
\label{F1.1.1}
\end{figure}

Then
the purpose of this paper is to study the disconnected manifold
$Z_{\mathbb R}$ of the set of Jacobi elements in ${\mathfrak g}$
associated to the generalized Toda lattices,
\begin{equation}
Z_{\mathbb R} =\left\{\ X=x+
  \displaystyle{{ \sum }_{ i=1}^l} (b_i e_{-\alpha _i } +  e_{\alpha _i})~ :~ x\in {\mathfrak h}, 
b_i \in {\mathbb R}\setminus \{0\},
 X\in S({\mathbb R}) \right\},
\nonumber
\end{equation}
its iso-spectral leaves $Z(\gamma)_{\mathbb R}$,
$\gamma\in \mathbb R^l$ and the construction of a smooth $connected$
compactification,  $ \hat Z(\gamma)_{\mathbb R} $ of each $Z(\gamma)_{\mathbb R}$.  The
construction of $\hat  Z (\gamma)_{\mathbb R}$  generalizes the  construction of
such a smooth compact manifold which  was carried out in \cite{KO-YE2} in the
important case of ${\mathfrak g}={\mathfrak{sl}}(l+1,{\mathbb R})$.
The construction there is based on the explicit solution structure
in terms of the so-called $\tau$-functions, which provide a local
coordinate system for the blow-up points. Then by tracing the solution orbit
of the indefinite Toda equation, the disconnected components in
$Z(\gamma)_{\mathbb R}$ are all
glued together to make a smooth compact manifold.  The result is maybe
well explained in Figure \ref{F1.1.1} for the case of ${\mathfrak{sl}}(3,{\mathbb R})$.
In the figure, the Toda orbits are shown as the dotted lines, and
each region labeled by the same signs in
$(\epsilon_1,\epsilon_2)$ with $\epsilon_i \in
\{\pm \}$ are glued together through the boundary (the wavy-lines) of
the hexagon. At a point of the boundary the Toda orbit
blows up in finite time, but the orbit can be uniquely
traced to the one in the next region (marked by the same letter
$A, B$ or $C$).
Then the compact smooth manifold $\hat Z(\gamma)_{\mathbb R}$
in this case is shown to be isomorphic to
the connected sum of two Klein bottles.
In the case of
${\mathfrak{sl}}(l+1, \mathbb R)$ for $l\ge 2$, $\hat Z(\gamma)_{\mathbb R}$ is shown to be
nonorientable and the symmetry group is the semi-direct product
of $({\mathbb Z}_2)^{l}$ and the Weyl group $W=S_{l+1}$, the permutation group.
One should also compare this
with the result in \cite{TO} where the compact manifolds are associated with
the definite (original) Toda lattice equation and the compactification
is done by adding only the subsystems. (Also see \cite{D-J-S} for some
topological aspects of the manifolds.)

The study of $Z(\gamma)_{\mathbb R} $  and of the   compact manifolds $ \hat
{Z}(\gamma)_{\mathbb R}  $  ,  can be  physically motivated by the appearance of the
indefinite Toda lattices in the context of symmetry reduction of the
Wess-Zumino-Novikov-Witten (WZNW) model. For example, the reduced system is
shown
in \cite{FE-TS} to contain  the indefinite Toda lattices. The compactification
${\hat Z}(\gamma)_{\mathbb R}$ can then be viewed as a concrete description of  (an
expected)
regularization of the integral manifolds of these indefinite Toda lattices,
where
infinities (i.e. blow up points) of the solutions of these Toda systems  glue
everything into  a smooth compact manifold.

In addition  our work is mathematically motivated by:
\begin{itemize}
\item[a)]
the work of Kostant in \cite{KOS} where he considered the real case with all
$b_i>0$,
\item[b)]
the construction of  the Toda lattice in \cite{LEZ-SAV}.
\end{itemize} 

In \cite{LEZ-SAV}, the solution $\{b_i(t) :i=1,\cdots,l\}$ in (\ref{1.1.2})
of the generalized Toda lattice equation
is shown to be expressed as an orbit on a connected
component $H_{\epsilon}$ labeled by
${\epsilon=(\epsilon_1,\cdots,\epsilon_l)}$ with $\epsilon_i=\pm 1$,
of the Cartan subgroup $H_{\mathbb R}$
defined in \S \ref{S2}, 
$$H_{\mathbb R}=\bigcup_{\epsilon\in {\{\pm 1\}^l}} H_{\epsilon}.$$
 This can be seen as follows:
 Let $g_{\epsilon}$ be an element of $H_{\epsilon}$ given by
\begin{equation}
g_{\epsilon}=h_{\epsilon}\exp\left(\displaystyle{\sum_{i=1}^l} f_ih_{\alpha_i}\right),
\label{1.1.7}
\end{equation}
which gives a map from $Z_{\mathbb R}$ into $H_{\mathbb R}$.
Here the element $h_{\epsilon}\in H_{\epsilon}$ satisfies
$\chi_{\alpha_{i}}(h_{\epsilon})=\epsilon_i$ for the group character
$\chi_{\phi}$ determined by $\phi \in \Delta$.
The solution $\{b_i(t):i=1,\cdots,l\}$ in (\ref{1.1.2}) 
is then directly connected to the group character
$\chi_{-\alpha_i}$ evaluated at $g_{\epsilon}$, i.e.
\begin{equation}
b_i(t)=\chi_{-\alpha_i}(g_{\epsilon}). 
 \nonumber
\end{equation}
The Toda lattice equation is now written as an evolution of $g_{\epsilon}$,
\begin{equation}
{d \over dt}g_{\epsilon}^{-1}{d \over dt}g_{\epsilon}=
\left[g_{\epsilon}^{-1}e_{+} g_{\epsilon}\ , \ e_{-}\right],  \nonumber
\end{equation}
where $e_{\pm}$ are fixed elements in the simple root spaces
${\mathfrak g}_{\pm\Pi}$ so that all the elements in ${\mathfrak g}_{\pm\Pi}$
can be generated by $e_{\pm}$, i.e.
${\mathfrak g}_{\pm\Pi}=\{Ad_h(e_{\pm}): h\in H\}$.  In particular,
we take
\begin{equation}
e_{\pm}=\displaystyle{\sum_{i=1}^l} e_{\pm\alpha_i}.
\nonumber
\end{equation}
Thus the Cartan subgroup $H_{{\mathbb R}}$ can be identified as
the position space (e.g. $f_i=q_i+\cdots+q_l$ for
${\mathfrak{sl}}(l+1,{\mathbb R})$-Toda lattice) of the generalized
Toda lattice equation, whose phase space is given by
the tangent bundle of $H_{\mathbb R}$.

One should also note that the boundary of each connected component
$H_{\epsilon}$
is given by either $b_i=0$ (corresponding to a subsystem) or $|b_i|=\infty$
(to a blow-up).
We are then led to our main construction of the compact smooth manifold
of $H_{\mathbb R}$ in attempting to generalize \cite{KOS} Theorem 2.4 to our
indefinite Toda case including $b_i<0$.

We define a set  $\hat H_{\mathbb R}$
containing the Cartan subgroup $H_{\mathbb R}$
by adding pieces corresponding to the blow-up points and the subsystems
(Definition  \ref{D5.1.1}).
Thus the set $\hat H_{\mathbb R}$ is defined as a disjoint union of
split Cartan subgroups $H_{\mathbb R}^A$ associated to semisimple portions
of Levi factors of all standard parabolic
subgroups determined by $A\subset \Pi$,
(Definition \ref{D2.2.2}),
$$
{\hat H_{\mathbb R}}=\displaystyle{\bigcup_{A\subset\Pi}\bigcup_{w\in W/W_{\Pi\setminus A}}
w\left(H_{\mathbb R}^A\right)\times\{[w]^A\}}$$
 The space
 $\hat H_{\mathbb R}$ then constitutes a kind of compact, connected
 completion of the disconnected Cartan
subgroup $H_{\mathbb R}$. Figure  \ref{F1.1.1} also describes how to connect the
connected components of $H_{\mathbb R}$ to produce
the connected manifold $\hat H_{\mathbb R}$ in the case of
$\mathfrak {sl}(3,{\mathbb R})$.
The signs must now be
interpreted as the signs of simple root characters
on the various connected components. The pair of signs $(+,+)$
corresponds to the connected component of
the identity $H$ and regions with a particular sign
represent one {\it single} connected component of $H_{\mathbb R}$. Boundaries
between regions with
a fixed sign  correspond to connected components
of Cartan subgroups arising from Levi factors of
parabolic subgroups. In addition the Weyl group $W$ acts on $\hat H_{\mathbb R}$.

Most of this paper is then devoted to describing the detailed structure of the
manifold $\hat H_{\mathbb R}$,
and in \S \ref{S8} we conclude with a diffeomorphism defined between ${\hat H}_{\mathbb R}$
and the isospectral manifold ${\hat Z}(\gamma)_{\mathbb R}$ as
identified with a toric variety $\overline{ (H_{\mathbb R}gB)}$
in the flag manifold ${\tilde G}/B$.

\subsection{Main Theorems\label{ss1.2}}

In connection with the construction ${\hat H}_{\mathbb R}$, we introduce in Definition
\ref{D3.1.1} the set of colored Dynkin diagrams.
 The colored Dynkin diagrams are simply Dynkin diagrams $D$ where some
of the vertices have been colored red or blue. For example in the case ${\mathfrak g}
={\mathfrak{sl}}(3,{\mathbb R})$:
$\circ_R-\circ$  (the sub-index $R$ indicates that $\circ$ is colored red).
The full set of colored Dynkin diagrams consists of pairs: $(D,[w]^{\Pi\setminus S})$ where
$D$ is a colored Dynkin diagram,   $S \subset \Pi$ denotes the set of
vertices that are colored in $D$ and $[w]^{\Pi\setminus S}$ is the coset of $w$ in
$W/W_S$. To each pair $(D,[w]^{\Pi\setminus S})$ one can
associate a {\it set} which is actually a {\it cell}. First
Definition \ref{D5.1.2} associates a subset of ${\hat H}_{\mathbb R}$ also denoted $
(D,[w]^{\Pi\setminus S})$.
\begin{figure}
\includegraphics{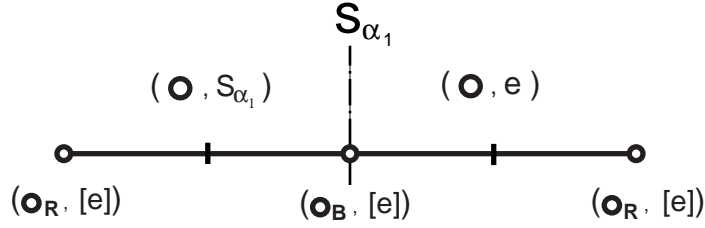}
\caption{The manifold $\hat H_{\mathbb R}$ parametrized by colored Dynkin diagrams
 for ${\mathfrak {sl}}(2, \mathbb R)$.}
\label{F1.2.1}
\end{figure}
This turns out to be a cell of codimension $k$
with $k=|S|$. We illustrate the example of ${\mathfrak{sl}}(2,
\mathbb R)$ in Figure \ref{F1.2.1}.

We consider the set ${\mathbb D}^k$ of colored Dynkin diagrams
$(D,[w]^{\Pi\setminus S})$ with $|S|=k$.  In Figure \ref{F1.2.1},
we have ${\mathbb D}^0=\{(\circ, e), (\circ, s_{\alpha_1})\},
{\mathbb D}^1=\{(\circ_B,[e]),(\circ_R,[e])\}$. We then obtain
the following theorem:

\begin{thm}
The collection of the sets $\{\ {\mathbb D}^k: k=0,1,\cdots,l \} $ 
gives a cell decomposition of $\hat {H}_{\mathbb R}$.
\label{T1.2.1}
\end{thm}

The chain complex ${\mathcal M}_*$ is introduced in \S \ref{S3}. The construction, in the case of ${\mathfrak{sl}}(2;{\mathbb R})$, is as follows:
\begin{equation}
{\mathcal M}_1={\mathbb Z}\left[ {\mathbb D}^0\right]
\overset{\partial_1}{\longrightarrow} {\mathcal M}_0 ={\mathbb Z}\left[{\mathbb D}^1 \right]
\nonumber
\end{equation}
Here the boundary map $\partial_1$
is given by 
\begin{equation}
\nonumber
\left\{
\begin{array}{ll}
\partial_1(\circ, e) &=(\circ_B, [e])-(\circ_R, [e]) , \\
\partial_1(\circ, s_{\alpha_1})  &=(\circ_B, [s_{\alpha_1}])-(\circ_R, [s_{\alpha_1}]).
\end{array}
\right.
\end{equation}
In particular, $(\circ):=\sum_{w\in W}(-1)^{\ell(w)}(\circ, w)=(\circ, e)-(\circ, s_{\alpha_1})$ is a cycle, where $\ell (w)$ is the
length of $w$, and $(\circ)$ represents the ${\hat H}_{\mathbb R}$.

The following theorem gives a topological description of  $\hat H_{\mathbb R}$:
\begin{thm}
The manifold $\hat {H}_{\mathbb R}$ is compact,
non-orientable (except if $\mathfrak g$ is of type $A_1$), and it has an action
of the Weyl
group $W$.  The integral
homology of $\hat {H}_{\mathbb R}$ can be computed as a ${\mathbb Z}[W]$ module as the
homology of the chain complex ${\mathcal M}_*$  in (\ref{3.3.6}).
\label{T1.2.2}
\end{thm}

Theorem  \ref{T1.2.2} is completed in Proposition \ref{P7.3.1}.  The $W$ action is (abstractly) introduced in Definition
\ref{D3.3.1} in terms of a representation-theoretic
induction process from smaller parabolic subgroups of $W$. The proof that
the $W$ action is well-defined is given in
Proposition \ref{P3.3.1} and Proposition \ref{P4.1.2}. The chain complex ${\mathcal M}^{CW}_*$ in Definition \ref{D7.3.1}
is defined so that it
computes integral homology. Since each $X_r\setminus X_{r-1}$ in Definition
\ref{D7.3.1} is a union of
cells $(D,[w]^{\Pi\setminus S})\in {\mathbb D}^{l-r}$, we
obtain an
identification between ${\mathcal M}^{CW}_*$ and the chain complex ${\mathcal M}*$.

Then using the Kostant map which can be described as a map from
$\hat H_{\mathbb R}$ into the flag manifold ${\tilde G}/B$ (in Definition \ref{D8.3.2}),
that is,
a torus imbedding, we obtain the following theorem:
\begin{thm}
The toric variety ${\hat Z}(\gamma)_{\mathbb R}$ is
a smooth
compact manifold which is diffeomorphic to ${\hat H}_{\mathbb R}$.
\label{T1.2.3}
\end{thm}
The complex version of this theorem has been proven in \cite{FL-HA2},
and our proof is essentially given in the same manner.

\subsection{Outline of the paper \label{outline}}

The paper is organized as follows:

In \S \ref{SSE5.2.2}, we begin with two fundamental examples, ${\mathfrak g}={\mathfrak{sl}}(l+1,{\mathbb R})$ for $l=1,2$, which summarize
the main results in the paper.

In \S \ref{S2}, we present the basic notations necessary for our discussions.
We then define a real group ${\tilde G}$ of rank $l$ whose split Cartan subgroup
 $H_{\mathbb R}$ contains $2^l$ connected components. We also define
 the Lie subgroups of $\tilde G$ corresponding to the subsystems
 and blow-ups of the generalized Toda lattice equations. The reason for
then introduction of $\tilde G$ and the Cartan subgroup $H_{\mathbb R}$ 
can be appreciated in Remark \ref{not-smooth} and Corollary \ref{nonsmooth}.

In \S \ref{S3}, we introduce colored Dynkin diagrams which will be shown
to parametrize the cells
in a cellular decomposition of the manifold ${\hat H}_{\mathbb R}$.
We then construct a chain complex ${\mathcal M}_*$ of the ${\mathbb Z}[W]$-modules
${\mathcal M}_{l-k}$ (Definition \ref{D3.3.2}). The parameters involved in
the statement of Theorem \ref{T1.2.1} are given here.

In \S \ref{S4} we show that Weyl group representations introduced in \S \ref{S3}
are well defined. In addition we define ${H}_{\mathbb R}^{\circ}$ in \S \ref{S4}
by adding some Cartan subgroups associated to semisimple
Levi factors of parabolic subgroups to $H_{\mathbb R}$.

In \S \ref{S5}, we define  $\hat H_{\mathbb R}$
as a union of the Cartan subgroup associated to semisimple Levi factors
of certain parabolic subgroups (Definition \ref{D5.1.1}) using translation
by Weyl group elements.
We also associate subsets of ${\hat H}_{\mathbb R}$ to colored
Dynkin diagrams. 

In \S \ref{S7},
we discuss the toplogical structure of ${\hat H}_{\mathbb R}$ expressing
${\hat H}_{\mathbb R}$ as the union of the subsets determined by the colored
Dynkin diagrams.
We then show that ${\hat H}_{\mathbb R}$ is a smooth compact manifold
and those subsets naturally determine a cell decomposition
(Theorem \ref{T1.2.2}).

In \S \ref{S8}, we consider a Kostant map between the isotropy subgroup $G^z_{\mathbb C}$
of $G_{\mathbb C}$ with $Ad(g)z=z$ and the isospectral manifold
$J(\gamma)_{\mathbb C}$ for some $\gamma \in {\mathbb R}^l$, which
can be also described as a map into the flag manifold ${\tilde G}/B$.
Then we show that the toric variety $\overline {(H_{\mathbb R} nB)}$
is a smooth manifold and obtain Theorem \ref{T1.2.3}.

\section{Examples of $\mathfrak {sl}(l+1, {\mathbb R})$, $l=1,2$}
  \label{SSE5.2.2}

This section contains two examples which are the source of insight for the main theorems in this paper. Most of the notation and constructions used later on in the paper can be anticipated  through these examples.

Our main object of study in this paper, ${\hat H}_{\mathbb R}$,  has nothing to do, in its construction, with the moment map and the Convexity Theorem of \cite{BFR:90}. However,  because of \cite{BFR:90} it  is expected to contain a polytope with vertices given by the action of the Weyl group which, gives rise to it through some gluings along its faces.  We will identify a convenient polytope of this kind in $\mathfrak h^{\prime}$ (dual of Cartan subalgebra) and show how all the pieces of ${\hat H}_{\mathbb R}$ would fit inside it.   The polytope in the example of ${\mathfrak{sl}}(3,{\mathbb R})$ is a hexagon. Since we are not dealing with the moment map in this paper, this is just  done for the purposes of motivation and illustration.  The polytope of the Convexity Theorem in \cite{BFR:90}, strictly speaking, is the convex hull of the orbit of $\rho$.  This sits inside the 
convex hull of the orbit of $2\rho$ which is all a part of ${\hat H}_{\mathbb R}$.

The  terms {\em dominant} and {\em antidominant} being relative to a choice of a Borel subalgebra, we refer to the chamber  in ${\mathfrak h}^{\prime}$ containing $\rho$ as {\em antidominant} for no other reason than the fact that we will make it correspond to  what we call the antidominant chamber of the Cartan subgroup $H_{\mathbb R}$ This odd convention is
 applicable throughout this section always  in connection with the  $2\rho$ polytope and our figures.

At this stage it is useful to keep in mind that:
\begin{enumerate}
\item  The bulk  -interior- of ${\hat H}_{\mathbb R}$ is made up of a split Cartan subgroup $H_{\mathbb R}$ which has $2^l$ connected components $H_{\epsilon}$ parametrized by signs $\epsilon=(\pm, \cdots , \pm)$ for ${\mathfrak{g}}$ with rank $l$, i.e. $
H_{\mathbb R}=\displaystyle{\bigcup_{\epsilon\in\{\pm 1\}^l}
H_{\epsilon}}.$ These disconnected pieces are glued together into a connected manifold by using Cartan subgroups  $H^A_{\mathbb R}$ associated to Levi factors of parabolic subgroups, which are determined by the set of simple roots $\Pi\setminus A$ for each $A\subset \Pi$. 
\label{E1}
\item The language of {\em colored Dynkin} diagrams and {\em signed colored} 
Dynkin diagrams is introduced in this paper in order to parametrize the pieces of 
${\hat H}_{\mathbb R}$. However the motivation for their introduction is that these 
diagrams parametrize the pieces of a convex polytope (hexagon) which lives in 
$\mathfrak h^{\prime}$, including its external faces and internal walls. The parameters for the pieces are colored and signed-colored Dynkin diagrams. 
\label{E2}
\item If one just looks at the antidominant  chamber intersected with the $2\rho$  polytope, then  it is easy to see that  it forms a {\em box}. This box is the closure of the  {\em antidominant chamber} $H_{\mathbb R}^{<}$ of the Cartan subgroup $H_{\mathbb R}$  (Definition \ref{D2.2.2}) inside $\hat H_{\mathbb R}$, i.e. 
${\hat H}_{\mathbb R}=\displaystyle{\bigcup_{w\in W}w\left( \overline {H_{\mathbb R}^{<}}\right) .}$
The box $\overline{H_{\mathbb R}^{<}}$ has {\em internal} chamber walls corresponding to simple roots (thus a Dynkin diagram appears) and it has {\em external} faces which are also parametrized by simple roots.  Hence one needs not just a Dynkin diagram but also two  {\em colors} to indicate  internal walls and external faces. The color {\em blue} indicates internal chamber walls and the color {\em red} indicates external faces.  {\em The $S$ is usually reserved for the set of colored vertices of a colored Dynkin diagram}.
\label{E3}
\item  To have a correspondence with the $2^l$ connected components of the Cartan subgroup $H_{\mathbb R}$ of diagonal matrices, the box  $\overline {H_{\mathbb R}^{<} }$ must be {\em further subdivided} into $2^l$ boxes $H_{\epsilon}^{<}$ with signs $\epsilon$, i.e. $
H_{\mathbb R}^{<}=\displaystyle{ \bigcup_{\epsilon\in{\{\pm 1\}^l}}H_{\epsilon}^{<}}. $
These boxes are parametrized by a Dynkin diagram where each simple root has a sign attached to it. The boundaries of these boxes are portions of the internal and external walls. Hence we use Dynkin diagrams with both signs and colors.  The boundaries between two signs require to be labeled as $0$. {\em The $A$ is usually reserved for  the set  of vertices of a  signed-colored Dynkin diagram assigned $0$ (subsystems)}.
\label{E4}
\item To {\em translate} the notation of  colored and signed-colored 
Dynkin diagrams to other chamber one needs to consider pairs $(D, w)$ where $w$ is a Weyl group element. However since the Weyl group action has non-trivial isotropy groups in portions of the polytope, it is necessary to consider Weyl group cosets $[w]^{\Pi\setminus S}$ for $S\subset \Pi$  the set of colored simple roots ( giving reflections generating an isotropy group). 
\label{E5}
\item To  {\em translate} the Levi subgroups around, note that the Weyl group of the Levi factor stabilizes the Cartan subgroup of the Levi factor
corresponding to the simple roots in $\Pi\setminus A$. Because of that we consider products $w(H_{\mathbb R}^A) \times \{ [w]^A \}\ $ so that $[w]^A$ is a coset in $W/W_{\Pi\setminus A}$,
$$
{\hat H}_{\mathbb R}=\displaystyle{\bigcup_{A\subset\Pi}\bigcup_{w\in W/W_{\Pi\setminus A}} w\left(H_{\mathbb R}^A\right)\times \{[w]^A\}.}
$$
\label{E6}
\end{enumerate}

\subsection{The example of $\mathfrak {sl}(2,{\mathbb R})$} \label{ss5.2}
The corresponding group in this example is given by ${\tilde G}=Ad(SL(2, {\mathbb R})^{\pm})$. The geometric picture
that corresponds to ${\hat H}_{\mathbb R}$ is a circle. Consider the interval
$[-2,2]$ where $-2$ and $2$ are identified.  Here $2$ represents $2\rho$.
Inside
this interval $(-2,2)$ we consider the subset $[-1,1]$. The points $-1,0,1$
divide  $[-2,2]$ into four open intervals. These open intervals will correspond
to the connected components of a Cartan subgroup $H_{\mathbb R}$ of $Ad(SL(2, {\mathbb R})^{\pm})$ when
the walls corresponding to the points 0 and 2 are deleted. Below we list each cell $w(H_{\epsilon}^{A,<})\times \{[w]^A\}$in ${\hat H}_{\mathbb R}$, 
where $H_{\epsilon}^{A,<}$ is the intersection of $H_{\epsilon}^A$
with the  strictly antidominant chamber (the superscript ${\leq}$ means that
the walls are included):

Let us take $h_{\alpha_1}={\rm diag}(1,-1)\in {\mathfrak h}$ and $h_{\epsilon}=Ad({\rm diag}(\epsilon,1))\in H_{\epsilon}$. Then any element
in $H_{\mathbb R}$ can be expressed as $h_{\epsilon}\exp(th_{\alpha_1})$ with
some parameter $t\in {\mathbb R}$. We denote $\exp(th_{\alpha_1})
={\rm diag}(a,a^{-1})$.

We first consider $A=\emptyset$. Then the  cell $H_+^{<}\times \{[e]\}$  is given by
\begin{equation}
\left\{ Ad({\rm diag}(a,a^{-1})) ~:~
 0<a<1 \right\}\times \{ [e]
\}=(\circ_+,e)  \leftrightarrow (0,1)
\nonumber
\end{equation}
Here $\chi_{\alpha_1}(h)=a^2$ for $h\in H_+$. Also we have the  set $H_-^{<}\times \{[e]\}$   as
\begin{equation}
\left\{ Ad({\rm diag}(-a,a^{-1})) ~:~ 0<a<1 \right\}\times
\{[e] \}=(\circ_-,e)  \leftrightarrow (1,2)
\nonumber
\end{equation}
with $\chi_{\alpha_1}(h)=-a^2$ for $h\in H_-$.

We now consider the case of $A=\{\alpha_1 \}=\Pi$, which
 corresponds to a subsystem of Toda lattice with $b_1=0$ in (\ref{1.1.5b}).
Then we have $H^{\Pi,\leq }_{\mathbb R}=\{ e \}$. This is the 
degenerate case of $A=\Pi$  which gives rise to the Levi factor of a Borel subgroup.  Since
the Levi factor does not contain a semisimple Lie subgroup,  $H^A_{\mathbb R}$
is defined to be  $\{ e \}$ (Definition \ref{D2.2.1}).
Here $[w]^{A}$  is just the element $w$. We have
\begin{equation}
\left\{ Ad(h_+)
\right\}\times\{e\}=(\circ_0,e) \leftrightarrow \{ 1 \}.
\nonumber  
\end{equation}
We now describe the box containing the strictly antidominant chamber
$H_{\mathbb R}^{<}$ of  $H_{\mathbb R}$. Since $H_{\mathbb R}^{<}$ is 
disconnected, $(\circ_0,e)$ has been used to glue the pieces together. We then have
a box given by
\begin{equation}
\begin{array} {ll}
& \displaystyle{(\circ, e)= (\circ_+,e)\cup (\circ_-,e)\cup  (\circ_0,e) 
\leftrightarrow (0,2).}
\end{array}
\label{5.2.4}
\end{equation}

The bijection which gives rise to local coordinates $\phi_e$ in Subsection 7.2
is given by either $\pm a^2$ or $0$. The set
$(\circ, e)$ is sent by $\phi_e$ to the interval $(-1,1)$.

We now apply $s_{\alpha_1}$ on $H_{\epsilon}^{<}\times\{[e]\}$ to obtain the cell in the $s_{\alpha_1}$-chamber
which corresponds to the
negative intervals:
\begin{equation}
\left\{ Ad({\rm diag}(a^{-1},a))~:~
 0<a<1 \right\}\times
\{ [e]
\}= (\circ_+, s_{\alpha_1})\leftrightarrow (-1,0).
\nonumber
 \end{equation}
with $\chi_{\alpha_1}(h)=a^{-2}$. However the local coordinate $\phi_{s_{\alpha_1}}$
is $\chi_{-\alpha_1}$ which equals $a^2$. We also have
\begin{equation}
\left\{ Ad({\rm diag}(a^{-1},-a))~:~
0<a<1 \right\}\times \{
[e] \}=(\circ_-, s_{\alpha_i}) \leftrightarrow (-2,-1).
 \nonumber
\end{equation}

Also for the case $A=\{ \alpha_1 \}$, and the set $s_{\alpha_1}(H^{A,\leq }_{\mathbb R})\times\{s_{\alpha_1}\}$  is given by
\begin{equation}
\left\{ Ad(h_+)
\right\}\times\{s_{\alpha_1}\}=(\circ_0,s_{\alpha_1} ) \leftrightarrow \{-1\}.
\nonumber
\end{equation}
Once again $\phi_{s_{\alpha_1}}$ is given as $0$, and we have
the set $s_{\alpha_1}(H_{\mathbb R}^{<})$, giving rise to an
open box
by gluing two disconnected pieces as before,
\begin{equation}
\begin{array} {ll}
& \displaystyle{ (\circ, s_{\alpha_1})= (\circ_+,s_{\alpha_1})\cup (\circ_-,s_{\alpha_1})\cup  (\circ_0,s_{\alpha_1})\leftrightarrow (-2,0).}
\end{array}
\nonumber
\end{equation}
\begin{figure}
\includegraphics{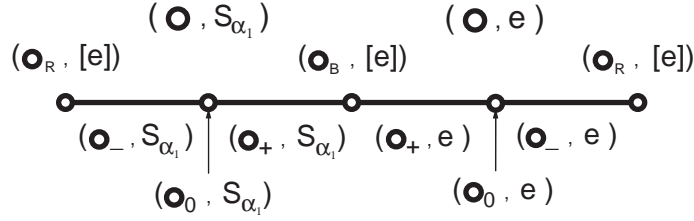}
\caption{The manifold $\hat H_{\mathbb R}$ parametrized by signed-colored Dynkin diagrams
 for ${\mathfrak {sl}}(2, {\mathbb R})$. The endpoints in the interval are identified giving rise to a circle.}
\label{F5.2.1}
\end{figure}

The image of $(\circ, s_{\alpha_1})$ under $\phi_{s_{\alpha_1}}$ is thus
$(-1,1)$. We also have the internal and external walls of the Cartan subgroup,
respectively:
\begin{equation}
\begin{array} {ll}
& \displaystyle{\left\{ Ad(h_+)
\right\}\times \{[e] \}=(\circ_B,[e]) \leftrightarrow\{ 0 \} }
\end{array}
\nonumber 
\end{equation}
\begin{equation}
\begin{array} {ll}
& \displaystyle{ \left\{ Ad(h_-)
\right\}\times \{ [e] \}= (\circ_R,[e])\leftrightarrow\{ 2 \}.}
\end{array}
\nonumber
\end{equation}
These are already associated to colored Dynkin diagrams $\circ_B$ and
$\circ_R$.

We can write down $\{ -2 \}$ by applying $s_{\alpha_1}$ as we did above.
However noting
$Ad({\rm diag}(-1,1))= Ad ({\rm diag}( 1,-1))$, we obtain the same set that
defines $\{ 2\}$. Thus $\{ 2 \}$ and $\{ -2 \}$ are identified.
We then obtain the
interval
$[-2,2]$ with $-2$ identified with $2$, which is $\hat H_{\mathbb R}$
diffeomorphic to a circle.  We illustrate this example in Figure \ref{F5.2.1}.

Maps can be easily found between the intervals on the left and the sets on the
right above.
With a suitable topology associated to the $(D,[w])$,topology defined
in terms of the coordinate functions $\phi_{e}$ and $\phi_{s_{\alpha_1}}$
in \ref{ss7.3},
 each interval or
point
on the left side is homeomorphic to the interval on the right. This is
what
is indicated with $  \leftrightarrow$.

\subsection{The example of  ${\mathfrak sl}(3, \mathbb R)$.}
\label{E5.2.2}

\begin{figure}
\includegraphics{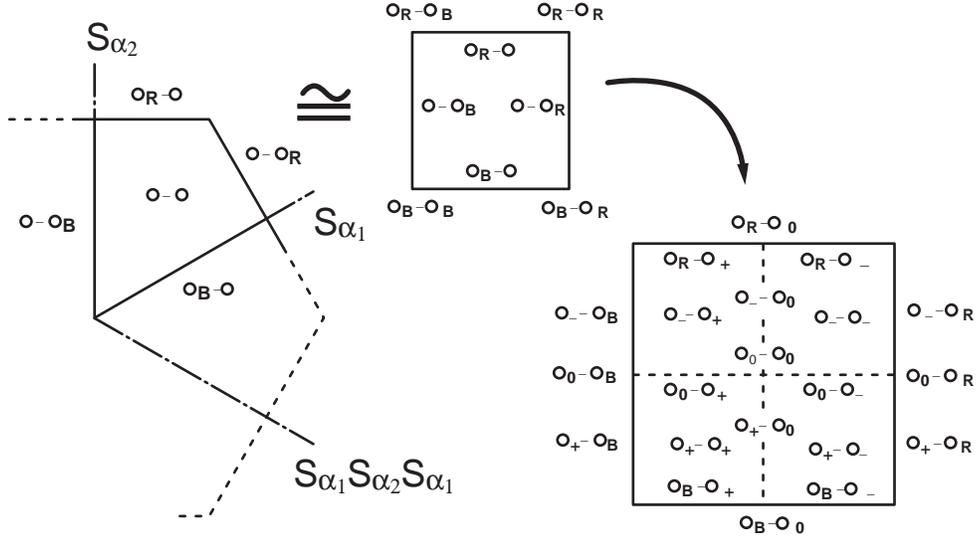}
\caption{ The {\it square} formed by the intersection of the antidominant chamber with the $2\rho$ hexagon for ${\mathfrak {sl}}(3, \mathbb R)$.}
\label{F5.2.2}
\end{figure}

We here  consider  the dual of its Cartan subalgebra ${\mathfrak
h}^\prime$ and, inside it,  a
convex region
bounded by a hexagon which is determined by the $W$ orbit of $2\rho$.  We will
later describe how to identify some of the faces of this hexagon. 
In Figure \ref{F5.2.2}, we illustrate the parametrization of the faces.
The two walls
of the antidominant chamber intersected with this convex region are denoted by
the colored Dynkin diagrams $\circ_B-\circ $ (the $s_{\alpha_1}$-wall) and
$\circ
-\circ_B$ (the $s_{\alpha_2}$-wall). The intersection of  two of the sides or
faces of the $2\rho$ hexagon with the antidominant chamber are  each
denoted by a
colored Dynkin diagram $\circ_R-\circ $ or $\circ -\circ_R$.  The four
colored Dynkin
diagrams $\circ_B-\circ $, $\circ -\circ_B$, $\circ_R-\circ $, $\circ -\circ_R$
form a {\it square}. This square is the intersection of the antidominant
chamber and the $2\rho$ hexagon. 

We now further subdivide this 
{\it square} into four {\it subsquares}
(see Figure \ref{F1.1.1} and Figure \ref{F5.2.2}).
Inside the $2\rho$ hexagon is the orbit of $\rho$ which gives rise to a new
smaller hexagon.  Consider the two faces in the smaller hexagon
intersecting the
antidominant chamber. Denote these (intersected with the antidominant
chamber) by
signed-colored diagrams. The $\circ _0 - \circ _+$ corresponds to the
unique face
intersecting the wall $\circ -\circ_B$, and $\circ _+ - \circ _0$
corresponds to
the unique face intersecting the wall $\circ_B-\circ $. Denote by $\circ
_0-\circ
_0$ the vertex of the $\rho$ hexagon which is the intersection of these two
faces. We now add a segment joining the vertex  $\circ _0-\circ _0$, to a point
in the interior of $\circ_R-\circ $, say the midpoint. Denote this segment
by the
signed-colored diagram $\circ _- - \circ _0$. We add another segment joining
$\circ _0-\circ _0$ to a point in the interior of $\circ -\circ_R$ and denote
this second segment by the signed-colored diagram $\circ _0-\circ _-$.  Now the
square is divided into four \lq\lq square \rq\rq regions denoted by $\circ
_{\pm}-\circ _{\pm}$.

Note that both $\circ_R-\circ $ and $\circ-\circ_R $ segments parametrized by
colored Dynkin diagrams, are now subdivided into two segments parametrized with
signed-colored Dynkin diagrams. For instance $\circ -\circ_R$ contains $\circ
_+-\circ_R$ and $\circ _--\circ_R$  and the intersection of these two is the
point $\circ _0-\circ_R$. Now the square $\circ _- - \circ _-$, for
example, has
a boundary which  consists of the segments parametrized by $\circ_R-\circ _-$,
$\circ _--\circ_R$, $\circ _0-\circ _-$ and $\circ _- -\circ _0$. The square
$\circ _--\circ _+$ has boundary $\circ_R-\circ _+$, $\circ _--\circ_B$, $\circ
_0-\circ _+$, $\circ _--\circ _0$. Here the $\alpha_2$-wall $\circ
-\circ_B$ has
also been subdivided into two pieces $\circ _+-\circ_B$ and $\circ
_--\circ_B$ by
the intersection with the $\rho$ hexagon. 

If we now consider the full set of signed-colored Dynkin diagrams by
translating
with $W$, we can fill the interior of the $2\rho$ hexagon with a total of $12$
regions. The four squares $(\circ _{\pm}-\circ _{\pm}, e)$ form the
intersection
of the antidominant chamber with the inside of the $2\rho$ hexagon.

\subsubsection{The sets in ${\hat H}_{\mathbb R}$ parametrized by the colored and signed-colored Dynkin diagrams}
 \label{hex2} 
We now proceed to describe explicitly some of the pieces of 
${\hat H}_{\mathbb R}$ corresponding to the signed-colored Dynkin diagrams 
$(D,[w]^{\Pi\setminus S})$ with $+,-$ or $0$ on the vertices in $\Pi\setminus S$ of the diagram 
$D$ and $[w]^{\Pi\setminus S} \in W/W_S$.

When $A=\emptyset$ implying no $0$'s in the vertices, we have
 $H^A_{\mathbb R}=H_{\mathbb R}$ which has $4$ connected
components,
\begin{equation}
{H}_{\mathbb R}=\displaystyle{\bigcup_{\epsilon\in \{\pm1\}^2}\left\{~h_{\epsilon}
{\rm diag}(a,b,c)~:~
a>0, ~b>0, ~abc=1~\right\}},
\label{hasl3}
\end{equation}
where $h_{\epsilon}=h_{(\epsilon_1,\epsilon_2)}={\rm diag}(\epsilon_2,\epsilon_1\epsilon_2,\epsilon_1)$ satisfying $\chi_{\alpha_i}(h_{\epsilon})=\epsilon_i$.
Since $A=\emptyset$ , $W_{\Pi\setminus A}=W$, and there is only one coset
$[e]^A=[e]$
in the ${\hat H}_{\mathbb R}$ construction. We now consider the signed-colored
Dynkin diagrams setting  $S=\emptyset$, that is, all the vertices are uncolored 
and they give a subdivision of the antidominant chamber inside the $2\rho$ hexagon.  
In order to move around this
chamber and its subdivisions using the Weyl group  we must  also consider six
elements $[w]^{\Pi} = w \in W$.
We then have, for $S=\emptyset$ and $A=\emptyset$, 
and $(\epsilon_1, \epsilon_2)=(\pm 1, \pm 1)$, 
\begin{equation}
\left(\circ _{\epsilon_1}-\circ _{\epsilon_2},e \right) 
=\displaystyle{
\left\{~h_{\epsilon}{\rm diag}(a,b,(ab)^{-1})~:~ab^{-1}<1,~ ab^2<1~\right\} 
\times \{[e]  \}}
\nonumber
\end{equation}
Note here that
the inequalities $|\chi _{\alpha_i}| < 1$  guarantee that the set  in question is contained in the chamber associated to $e$.
Here the local coordinate $\phi_e$ is given by $(\chi_{\alpha_1},\chi_{\alpha_2})$ which equals
$(ab^{-1},ab^2)$.

\medskip

For $A=\{\alpha_1\}$, we have $h_{\epsilon}={\rm diag}(\epsilon_2,\epsilon_2, 1)\in H_{\epsilon}^A$ with $\epsilon=(1,\epsilon_2)$; and we  multiply this element  with 
$ {\rm diag} (1, a,a^{-1})$. This is a typical element in the connected component of the Cartan associated to the Levi factor,  in accordance with the definition of $H^A_{\mathbb R}$ in Definition \ref{D2.2.2}. This gives:
$$(\circ _0-\circ _{\epsilon_2}, e)=\left\{\ {\rm diag} (\epsilon_2, \epsilon_2 a,a^{-1}) 
\ : \
 0<a < 1 ~\right\}\times \{ [e]^A  \} $$
where the local coordinate $\phi_e$ is given by $(0,\epsilon_2 a^2)$.

\medskip

For $A=\{\alpha_2\}$, we have in a similar way with $h_{(\epsilon_1,1)}={\rm diag}(1, \epsilon_1, \epsilon_1)$:
$$(\circ _{ \epsilon_1}-\circ _0, e)=\left\{\ {\rm diag }(a, \epsilon_1a^{-1},\epsilon_1) \ : \
 a > 0 , |a| < 1 \right\}\times \{[e]^A  \}$$
where $\phi_e$ now equals $(\epsilon_1 a^2,0)$.

\medskip

We have an open  {\em square} associated with the interior of the antidominant chamber,
\begin{equation}
\begin{array} {ll}
&  \displaystyle{ }\\
& \displaystyle{(\circ-\circ,e)= \underset {(\nu_1,\nu_2)\in \{\pm 1,0 \}^2
} {\bigcup }~(\circ_{\nu_1}-\circ_{\nu_2}, e)}~.
\end{array}
\label{5.2.20}
\end{equation}
The image of this set under the map $\phi_e$ is an open square $(-1,1)\times (-1,1)$.

We now write down the boundary of this
 square: We here give an explicit form of $(\circ_R-\circ,[e]^{\Pi \setminus S})$, and the others can be obtained in the similar way.
We first have, for $S=\{ \alpha_1 \}, A=\emptyset $, so that
$[e]^A=[e]=[s_{\alpha_i}]$ for $i=1,2$,
\begin{equation}
(\circ_R-\circ _{\epsilon_2}, [e]^{\{\alpha_2\}})=\left\{\ {\rm diag} \left( \epsilon_2 a, -\epsilon_2 a,-a^{-2}\right) \ :\
 0<a<1 \right\}\times \{ [e]  \}
\nonumber
\end{equation}
Here $\phi_e$ equals $(\chi_{\alpha_1},\chi_{\alpha_2})$, and is given by
$(-1, \epsilon_2 a^3)$.

With $A= \left\{ \alpha_2 \right\}$, we have:
\begin{equation}
(\circ_R-\circ_0,[e]^{\{\alpha_2\}})=\{\ {\rm diag} (1,-1,-1)~ \} \times \{ [e]^{\{\alpha_2\}}
\}
\nonumber
\end{equation}
The map $\phi_e$ is $(\chi^{\Delta^A}_{\alpha_1},0)$ and equals $(-1,0)$.
We then have
\begin{equation}
 \displaystyle{(\circ_R -\circ, [e]^{\{\alpha_2\}})=
 \bigcup_{\nu\in\{\pm 1,0\}}(\circ_R-\circ _{\nu}, [e]^{\{\alpha_2\}})}~.
\nonumber
\end{equation}
The image of this set under $\phi_e$ is thus $\{ -1 \}\times (-1,1)$.

\bigskip

We now consider the parts of the Cartan subgroup of Levi factors corresponding to other chambers inside the hexagon.
Note that if we apply $s_{\alpha_1}=w$ to $(\circ_{\epsilon_1}-\circ_{\epsilon_2},e)$, we obtain
\begin{equation}
\displaystyle{ \left(\circ_{\epsilon_1}-\circ_{ \epsilon_1\epsilon_2},s_{\alpha_1}\right)=\left\{ {\rm diag} \left(\epsilon_1\epsilon_2 b, \epsilon_2 a,\epsilon_1 (ab)^{-1}\right)  :
 ab^{-1} <1,  ab^2< 1  \right\}\times \{ [e]  \}. }
\nonumber
\end{equation}
Since $sign (\chi_{\alpha_1})= \epsilon_1$ and $sign (\chi_{\alpha_2})=\epsilon_1\epsilon_2$, this 
set is no longer contained in $H_{(\epsilon_1, \epsilon_2)}$ but rather in $H_{(\epsilon_1,\epsilon_1 \epsilon_2)}$. This justifies the notation $
(\circ _{\epsilon_1}-\circ _{\epsilon_1\epsilon_2})$. One should note that $\epsilon$ for the component $H_{\epsilon}$ in $H_{\mathbb R}$
did not  change when one uses the simple roots associated to the {\em new} positive system
$s_{\alpha_1}\Delta_+$.

Also notice that for $S=\{\alpha_1\}$ we have
\begin{equation}
\left(\circ_R-\circ _+,[s_{\alpha_1}]^{\Pi\setminus S}\right) =\left\{\ {\rm diag} \left(-a,
a,{-a^{-2}}\right)\ : \ a>0 \right\}\times \{[e]  \}
\nonumber
\end{equation}
and we have an identification of $(\circ_R-\circ _+,[s_{\alpha_1}]^{\Pi\setminus S})$ and $ (\circ_R-\circ _-, [e]^{\Pi\setminus S})$. They are the same set. In fact, now $S=\{\
\alpha_1 \}$ and $e$ and $s_{\alpha_1}$ give the same coset in $W/W_S$ so the
corresponding signed-colored Dynkin diagrams agree too.
Similarly $ (\circ_R-\circ _-,[s_{\alpha_1}]^{\Pi\setminus S})$ is the same as
$(\circ_R-\circ _+, [e]^{\Pi\setminus S})$. In our hexagon this means that {\em two of the outer walls must be glued}. The fact that a segment with a sign $+$ is glued to one with $-$ corresponds to the fact that the two contiguous chambers will form a {\em Mobius band} after the gluing {see Figure \ref{FMoeb}).
\begin{figure}
\includegraphics{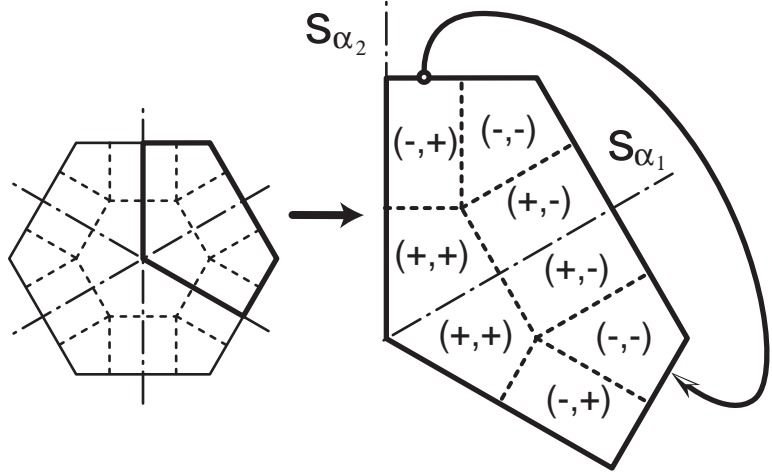}
\caption{Gluing creates a Mobius band out of  two contiguous chambers }
\label{FMoeb}
\end{figure}
What this means is that in our  geometric picture consisting of the inside
of the
$2\rho$ hexagon, some portions of the boundary need to be identified. Such
identifications take place on all the chambers.  This identification
provides the gluing rule given in Lemma 4.2 in \cite{KO-YE2}
for the case of ${\mathfrak{sl}}(n, {\mathbb R})$.

\subsection{The chain complex ${\mathcal M}_*$}

 We describe  ${\mathcal M}_*$ in 
terms of colored Dynkin diagrams. The $\mathbb Z$ modules of chains ${\mathcal M}_k$ are then given by
(see Figure \ref{F3.3.1}):
\begin{itemize}
\item 
${\mathcal M}_2={\mathbb Z}\left[ (\circ-\circ, w):
w\in W\right]$. The cells are the {\em chambers of the 
$2\rho$ hexagon}, and ${\rm dim}{\mathcal M}_2=6$.

\item   ${\mathcal M}_1$ consists of the cells 
parametrized by the
colored
Dynkin diagrams $(\circ_R-\circ, [w]^{\{\ \alpha_2 \}})$, $(\circ_B-\circ,
[w]^{\{\ \alpha_2 \}})$ with $w\in\{e,
 s_{\alpha_2}, s_{\alpha_1}s_{\alpha_2}\}$
and $(\circ-\circ_R,[w]^{\{\ \alpha_1 \}})$, $(\circ-\circ_B,[w]^{\{\
\alpha_1 \}})$ with
$w\in\{e, s_{\alpha_1}, s_{\alpha_2}s_{\alpha_1}\}$.  These are the {\em sides  of the different chambers
of the hexagon}.  The {\em blue} ($B$) stands for {\em internal chamber wall} and
the {\em red} ($R$) for {\em external face of the hexagon}. The dimension is then given by ${\rm dim}{\mathcal M}_1=12$.

\item  ${\mathcal M}_0={\mathbb Z}\left[(\circ_s-\circ_t, [e]): s,t\in\{R,B\}\right]$.
These are the four vertices of a chamber. Because of identifications, 
there are only four different points coming from all the six chambers, that is, ${\rm dim}
{\mathcal M}_0=4$.
\end{itemize}

\begin{figure}
\includegraphics{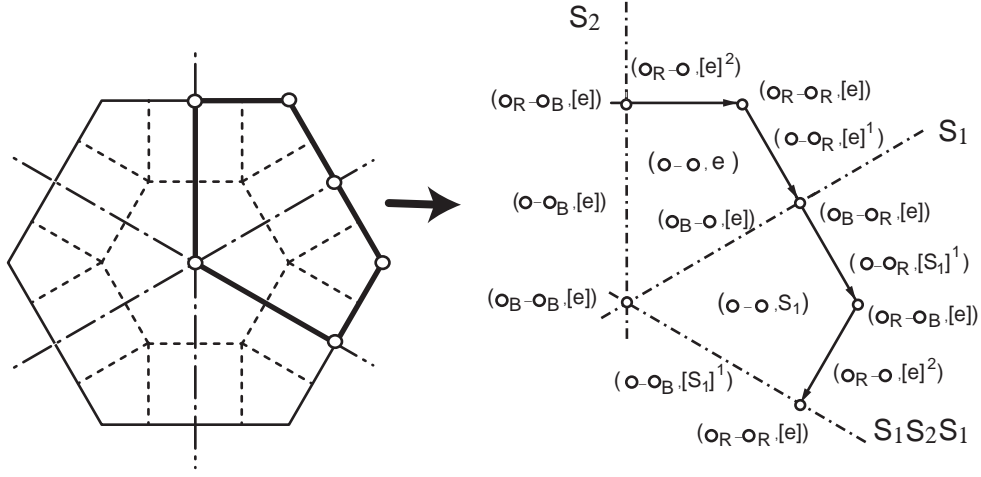}
\caption{The manifold $\hat H_{\mathbb R}$ parametrized by colored Dynkin diagrams
 for ${\mathfrak {sl}}(3, \mathbb R)$.}
\label{F3.3.1}
\end{figure}

The chain complex ${\mathcal M}_* : {\mathcal M}_2
\overset {\partial_2}{\longrightarrow}{\mathcal M}_1
\overset {\partial _1}{\longrightarrow}{\mathcal M}_0$
leads to the following integral homology: $H_2=0$ , $H_1={\mathbb Z}^3\oplus{\mathbb
Z} /2 {\mathbb Z}$, $H_0={\mathbb Z}$.  This implies that
$\hat H_{\mathbb R}$ is nonorientable and is equivalent to the
connected sum of two Klein bottles. Also note that the Euler character is
$6-12+4=-2$.
 According to Proposition \ref{P7.3.1} this computes the homology
$H_*(\hat H_{\mathbb R}, {\mathbb Z})$
of the compact smooth manifold ${\hat H}_{\mathbb R}$.
The torsion ${\mathbb Z}/2 {\mathbb Z}$ in $H_1$ has the following representative:
\begin{equation}
\begin{array}{ll}
c_1=&\displaystyle{\sum_{w\in W/W_{\{s_{\alpha_2}\}}}
(-1)^{\ell(w)}\left(\circ_R-\circ,[w]^{\{s_{\alpha_2}\}}\right)} \\
& \quad - \displaystyle{\sum_{w\in W/W_{\{s_{\alpha_1}\}}}
(-1)^{\ell(w)}\left(\circ_R-\circ,[w]^{\{s_{\alpha_1}\}}\right)}~.
\end{array}
\nonumber
\end{equation}
Here $\ell(w)$ denotes the length of $w$. If we let
\begin{equation}
c_2=\underset { w\in W} {\sum} (-1)^{\ell(w)} (\circ-\circ, w),
\nonumber
\end{equation}
then $\partial_2 (c_2)=2c_1$. 

From the chain complex one can compute  the three dimensional vector space
$H^1( {\hat H}_{\mathbb R}, {\mathbb Q} )$ as a $W$ module. This is a direct sum
of a one dimensional non-trivial (sign) representation and
the two dimensional reflection representation. The representation
$H_0( {\hat H}_{\mathbb R}, {\mathbb Q})$ is trivial.

\section{Notation, the group $\tilde G$}
 \label{S2}

\subsection {Basic Notation} \label{ss2.1}

The important Lie group for the purposes of this  paper is a group $\tilde G$
that will be technically introduced in subsection \ref{ss2.2}. This group is $\mathbb R$ split
and has
a split Cartan subgroup with $2^l$ components with $l={\rm rank}(G)$. For example $SL(3,\mathbb R)$
is of this
kind but $SL(4,\mathbb R)$ is not. We need to
introduce the following standard objects before it is possible to define and
study $\tilde G$.
\begin{notation} {\it Standard Lie theoretic notation:}
\label{D2.1.1}
We adhere to standard notation and to the following conventions: ${\cdots} _{\mathbb C}$ denotes a complexification; for any set $S\subset \Pi$, ${\cdots }_{S}$ is an object related to a parabolic subgroup or subalgebra determined by $S$; for any subset $A\subset \Pi$ , ${\cdots}^A$ is associated to the parabolic determined by $\Pi\setminus A$. However we will not simultaneously employ ${\cdots }^A$ and ${\cdots}_S$  for the same object ${\cdots}$.
As in the introduction  $\mathfrak g$ denotes a
real split semisimple Lie algebra of rank $l$ with complexification$\mathfrak g_{\mathbb C}={\mathfrak g} \underset {\mathbb R} {\otimes }
{\mathbb C}$. We also have   $G_{\mathbb C}$ the connected adjoint
Lie group with Lie algebra $\mathfrak g_{\mathbb C}$ ($G_{\mathbb C} = Ad(G^s_{\mathbb C})$) and $G$ the connected real semisimple Lie subgroup of  $G_{\mathbb C}$ with Lie algebra $\mathfrak g$. Denote $G^s_{\mathbb C}$  the simply connected complex Lie group associated to ${\mathfrak g}_{\mathbb C}$. We
list some additional very standard Lie theoretic notation:
\begin{itemize}
\item  ${\mathfrak g}^{\prime}={\rm Hom}_{\mathbb R} ({\mathfrak g}, {\mathbb
R})$ and ${\mathfrak g}^\prime_{\mathbb C}={\rm Hom}_{\mathbb C} ({\mathfrak g}_{\mathbb C},{\mathbb C})$.

\item Given $\lambda \in {\mathfrak g}^{\prime}$ and $x\in {\mathfrak g}$, 
$\langle \lambda, x\rangle $ is $\lambda$ evaluated in $x$, $\langle \lambda, x\rangle=\lambda(x)$. 

 \item $( , )$ the bilinear form on $\mathfrak g$ or ${\mathfrak g}_{\mathbb C}$  given by
the Killing form (the same notation applies to the Killing form on $\mathfrak
g^\prime$ and $\mathfrak g_{\mathbb C}^{\prime}$).

\item  $\theta$ a Cartan
involution on $\mathfrak g$.
 
\item  $ {\mathfrak g}={\mathfrak k}+{\mathfrak p}$ the Cartan
decomposition of $\mathfrak g$ associated to $\theta$, where
$\mathfrak k$ is the Lie algebra of a  maximal compact
subgroup $K$ of the adjoint group $G$ and $\mathfrak p$ is the orthogonal
complement to $\mathfrak k$ with respect to the Killing form. 

\item  $e_\phi\in {\mathfrak g}$, $\mathfrak h$ root vectors chosen so
that $(e_{\phi},e_{-\phi})=1$.

\item  $\Delta_+ \subset \Delta$ be a
fixed
system of positive roots.

\item  ${\mathfrak b}= {\mathfrak h} +  \underset {
{\phi} \in \Delta_+ }  { \sum  }{\mathbb R } e_{\phi} \,$ (Borel subalgebra).

\item  ${\mathfrak n}=[{\mathfrak b},{ \mathfrak b}] =\underset { \phi \in \Delta_+ } 
{ \sum }{ \mathbb R} e_\phi \, $ and $\bar {\mathfrak n}=\underset { \phi \in \Delta_+ }   {
\sum  } {\mathbb R} e_{-\phi }\, $.

\item  ${\mathfrak n}_{\mathbb C}={\mathfrak n}
\underset {\mathbb R}{\otimes }{\mathbb C}$, $\bar {\mathfrak n}_{\mathbb C}=\bar {\mathfrak n }\underset {\mathbb R}{\otimes } {\mathbb C}$.
\item  $H_{\mathbb C}$ Cartan subgroup of $G_{\mathbb C}$
with Lie algebra ${\mathfrak h}_{\mathbb C}$.

\item  $H=H(\Delta)$ the connected Lie subgroup of $G$
with Lie algebra ${\mathfrak h}$, $H=\exp({\mathfrak h})$.

\item  $H^1_{\mathbb R}$ Cartan subgroup of $G$ with Lie algebra $\mathfrak h$.

\item   $W$ the Weyl group of $\mathfrak g$ with respect to $\mathfrak h$ or
the Weyl group of $H^1_{\mathbb R}$.

\item   $W_S\subset W$, the group
generated by the simple reflections corresponding to the elements in $S\subset \Pi$.
\end{itemize}
\end{notation}
\begin{rem}
The group $H^1_{\mathbb R}$ \cite{WA}  p.59, 2.3.6  consists of all $g\in G$ such that
$Ad(g)$ restricted to $\mathfrak h$ is the identity.  This Cartan subgroup will be
usually disconnected and $H$ is the connected component of the identity $e$.
The Weyl group $W$ of the Cartan subgroup $H^1_{\mathbb R}$ is isomorphic to the
group
which is generated by the simple reflections $s_{\alpha_i}$ with $\alpha_i \in
\Pi$ and which agrees with the Weyl group associated to the pair $({\mathfrak g}_{\mathbb
C}, {\mathfrak h}_{\mathbb C})$. This is because our group is assumed to be $\mathbb R$ split.
\end{rem}

\begin{example}
\label{E2.1.1}
The reader may wish to read the whole paper with the following well-known
example
in mind. Let $G_{\mathbb C}^s= SL(n,\mathbb C)$ and  $G_{\mathbb C}= Ad (SL(n,{\mathbb
C}))$. The
second group is obtained by dividing $SL(n,{\mathbb C})$  by the finite abelian
group
consisting of the $n$ roots of unity times the identity matrix.   We can set
$G=Ad (SL(n,{\mathbb R}))$. Note that if $n$ is odd then $SL(n, {\mathbb R})=Ad
(SL(n, {\mathbb R}))=G$. If $n$ is even then $Ad (SL(n,{\mathbb R} ))$ is obtained by
dividing by $\{ \pm I \}$ ($I$ the identity matrix). In this  example, $\mathfrak g$ consists of
traceless $n\times n$ real matrices.
The Cartan subalgebra $\mathfrak h$ in the $Ad (SL(n, {\mathbb R}))$ case can be taken to
be the space of traceless $n\times n$ real diagonal matrices. The root vectors
$e_\phi$ are the various matrices with all entries $a_{i,j}=0$ if $(i,j)\not
=(i_o,j_o)$ and $a_{i_o,j_o}=1 $ where $i_o\not=j_o$ are fixed integers.
 In this case of $ G= Ad(SL(n,
{\mathbb R}))$, with $\mathfrak h$ chosen as above, the Cartan subgroup $H^1_{\mathbb R}$ of
$G$ consists of $Ad$ applied to the group of  real diagonal matrices of
determinant one. The group $H$ consists of $Ad$ applied to all diagonal
matrices
${\rm diag}(r_1,..,r_n)$ with $r_i > 0$ the connected component of the identity of
$H^1_{\mathbb R}$.
\end{example}
\begin{notation}
\label{D2.1.2}
The following elements in ${\mathfrak h}$ and its dual ${\mathfrak h}'$
will appear often in this paper:
\begin{itemize}
\item  $\check \alpha_i= \frac {2 \alpha_i}{(\alpha_i , \alpha_i)}$
the coroots.
\item  $m_{\alpha_i}$ $i=1,..,l$, the fundamental weights,
$(m_{\alpha_i}, \check \alpha_j)= \delta_{i,j}$.
\item  $ m_{\alpha_i}^\circ $ the unique
element in $\mathfrak h$ defined by $\langle m_{\alpha_i}, x\rangle
=(m_{\alpha_i}^\circ , x)$.
\item  $h_{\alpha_i}$ is the unique element in $\mathfrak h$ such that
$(h_{\alpha_i},  x)=\langle \check \alpha_i, x\rangle $.
\item  $y_i=\frac {2\pi m_{\alpha_i}^\circ}{(\alpha_i,
\alpha_i)}$.
\item  ${\mathfrak h}^{<}=\{\ x\in {\mathfrak h}: \langle \alpha_i,x\rangle  < 0
\text{
for all } \alpha_i\in \Pi \} $.
\end{itemize}
\end{notation}
\begin{example}
\label{E2.1.2}
 Consider  $G_{\mathbb C}= Ad(SL ( n, {\mathbb C}))
\subset Ad(GL(n, {\mathbb C}))$ and $n$ even . We view ${\mathfrak g}_{\mathbb C}$ inside
the Lie algebra of $GL(n,{\mathbb C})$ (as $n\times n$ complex traceless matrices).
Then $m_{\alpha_i}^\circ=
{\rm diag}(t_1,..,t_n)+z$ where $t_j=1$ for $j\leq i$ and $t_j=0$ for $j>i$. The
element $z$ is in the center of the Lie algebra of $GL(n, {\mathbb C})$ and thus
$ad(z)=0$.  For instance if $n=2$ then $m_{\alpha_1}^\circ ={\rm diag} (\frac
{1}{2},-\frac {1}{2})= {\rm diag }(1,0)-{\rm diag}(\frac {1}{2},\frac {1}{2})$.
\end{example}

\subsection {The group $\tilde G$}
 \label{ss2.2} 

We now define, following \cite{KOS}, 3.4.4 in p.241 an enlargement $\tilde G$
of the split group $G$. The purpose of this is to \lq\lq complete\rq\rq the
Cartan subgroup $H^1_{\mathbb R}$ forcing it to have $2^l$ connected components
where
$l$ is the rank of $G$. In the case $G=SL(l+1,{\mathbb R})$ where $l$ is even, $\tilde
G=G$ already.

We thus  consider another real split group, slightly bigger than $G$. Let $x,y
\in \mathfrak h$ so that $x+\sqrt {-1} y \in \mathfrak h_{\mathbb C}$. We recall the conjugate
linear automorphism of ${\mathfrak g}_{\mathbb C}$ given by $z^c= x-\sqrt {-1} y$ where $z=x+\sqrt {-1} y$ and
$x,y\in { \mathfrak g}$. We recall that this automorphism induces an automorphism
$G_{\mathbb C} \to G_{\mathbb C}$ , $g\mapsto g^c$. We  thus let
\begin{equation}
{\tilde G}=\{\ g\in G_{\mathbb C}: g^c=g \}.
\nonumber
\end{equation}
By  \cite{KOS} Proposition 3.4, we also have:
\begin{equation}
{\tilde G}=\{\ g\in G_{\mathbb C}:Ad(g){\mathfrak g} \subset  {\mathfrak g} \}.
\nonumber
\end{equation}
Then we have the following Proposition whose proof will be given later (right after Proposition
\ref{P2.2.2}):
\begin{prop}
Let $G= Ad(SL(n,\mathbb R))$ ,then $\tilde G
\cong
SL(n,\mathbb R)$ for $n$ odd and $\tilde G \cong Ad (SL(n,\mathbb R)^{\pm}   )$
if $n $ is
even.  Thus $\tilde G$ is disconnected whenever $n$ is even.\label{P2.2.1}
\end{prop}

We now describe an element $h_i$ in the Cartan subgroup of $\tilde G$. First we have:
\begin{lem}
Let $x,y\in \mathfrak h$ and $\alpha_i\in \Pi$. Then  the
numbers $\exp( \langle \alpha_i, x+ \sqrt {-1} y\rangle    )$ are  real  if and only
if $y$ has the form:
$y=  { \sum_{i=1}^l  } k_iy_i \,$ where $k_i$ is an integer
and $y_i$ is as in  Notation \ref{D2.1.2}.  The elements  $h_i=\exp(\sqrt {-1} k_iy_i )$ with
$k_i$ odd
are in $\tilde G$ and satisfy $h_i^2=e$ with $h_i\not=e$.\label{L2.2.1}
\end{lem}

\begin{Proof}
 It is enough to consider the case when $x=0$. Suppose
first that $y=  { \sum_{i=1}^l  } k_i y_i\,$
where each $k_i$ is an integer.
Then $e^{\sqrt {-1} \langle \alpha_i,  y\rangle}=e^{ \sqrt {-1} \pi
k_i}$ takes either $+1$ or $-1$. Conversely, suppose that  all the $e^{\langle \alpha_i,
\sqrt {-1} y\rangle }$ with $i=1,\cdots, l$ are
real.
Then $e^{\langle \alpha_i, \sqrt {-1} y\rangle }$ equals $\pm 1$ and $\langle
\alpha_i , y\rangle = k_i \pi$ for all
$i=1,..,l$. This implies that , for each $i$, $\langle \check \alpha_i,
y\rangle =\frac
{2 k_i
\pi }{(\alpha_i, \alpha_i) }$. Since $\mathfrak g$ is semisimple and the
$\{m_{\alpha_i}^\circ :i=1,..,l\}$ forms a basis of $\mathfrak h$,  necessarily $y= { \sum_{i=1}^l  } c_i m_{\alpha_i}^\circ \,$ with
$\langle \check
\alpha_i,  y\rangle =c_i$ and $\langle \check \alpha_i ,y\rangle = \frac
{2k_i \pi}{(\alpha_i,
\alpha_i) }$.  This proves the first part of the statement in Lemma
\ref{L2.2.1}.  We
note that  since  all the $e^{\langle \alpha_j,y_i\rangle }$ are real  then
$e^{\langle  \phi, \sqrt {-1} y_i\rangle}$ is also real for any root
 $\phi$ which is not necessarily simple.
Therefore each $Ad(h_i)$
stabilizes all the root spaces of $\mathfrak h$  in $\mathfrak g$. Since $\mathfrak h$
is also
stabilized, we obtain that $Ad(h_i) ({\mathfrak g}) \subset  {\mathfrak g}$ and thus $h_i\in
{\tilde G}$. Clearly $h_i= \exp(\sqrt {-1} y_i)$ satisfies $h_i^2=1$ where $h_i\not=e$.
This
is because $h_i$ is representable by a diagonal matrix with entries of the form
$\pm 1$. Moreover, at least one diagonal entry must be equal to $-1$.
\end{Proof}

\begin{example}
In the case of $Ad(SL(n, {\mathbb R})^{\pm})$
with $n$ even, $h_i=\exp(\sqrt {-1} y_i)$ is just $Ad ({\rm diag} (r_1,..,r_n))$
where $r_j=1$ if $j \leq i$ and $r_j=-1$ if $j > i$.
When $n$ is odd then we have $h_i=(-1)^{n-i} {\rm diag} (r_1,..,r_n)$
with the same notation as above for the $r_i$.
In this case $SL(n,{\mathbb R})=Ad( SL(n, {\mathbb R}))$.
\label{E2.2.1}
\end{example}

We now describe a split Cartan subgroup $H_{\mathbb R}$
 of $\tilde G$ and other items related to it. The $H_{\mathbb R}$ is
 the real part of $H_{\mathbb C}$ on $\tilde G$,
\begin{equation}
  H_{\mathbb R}=H_{\mathbb C} \cap {\tilde G},
\nonumber
\end{equation}
which has $2^l$
components (see Proposition \ref{P2.2.2} below). We denote  by $B$  a Borel subgroup with Lie algebra $\mathfrak h +\mathfrak n$ contained in $\tilde G$ or $G$ as will be clear from the context. Thus in $\tilde G$ this is $H_{\mathbb R}N$, $N=\exp (\mathfrak n)$.  From the Bruhat decomposition applied to $\tilde G$, we have
\begin{equation}
{\tilde G}=\bigcup_{w\in W}\bar N \hat w H_{\mathbb R}N ,
\nonumber
\end{equation} 
with $\bar N=\exp (\bar {\mathfrak n})$. Here
$\hat w$   stands for any representative of the Weyl group 
element $w\in W$ in the normalizer of the Cartan subgroup.   Keeping this
in mind we will harmlessly drop the $\hat{}$ from the notation. 

In addition let $\chi_{\phi}$ denote the group character determined by $\phi\in \Delta$;
on $H_{\mathbb R}$ each $\chi_{\phi}$ is real and cannot take the value zero.
Thus a group character has a fixed sign on each connected component. We denote by $sign(\chi_{\phi}(g))$ the sign of this character on a specified element $g$ of $H_{\mathbb R}$. We let  ${\mathcal E}= {\mathcal E}(\Delta)=\{\pm 1\}^l=\{\
(\epsilon_1,..,\epsilon_l) : \epsilon_i\in \{ \pm 1 \}, i=1,..,l\}$.
Then the set $\mathcal E$ of all $2^l$ elements $\epsilon$ (or functions $\epsilon
:\Pi\to
\{\pm  1\}$) will parametrize connected components of the Cartan subgroup
$H_{\mathbb R}$ of $\tilde G$ and corresponds to the signs $\epsilon_k= s_k s_{k+1}$
in the indefinite Toda lattice  in p.323 of  \cite{KO-YE1}.
In connection with the connected components let  $h_\epsilon = \prod_{\epsilon_i\not=1} h_i$ for $\epsilon \in
{\mathcal E}$ and $h_i=\exp(\sqrt {-1} y_i)$ in Lemma \ref{L2.2.1}.
Now  $H_\epsilon=h_\epsilon H =\{\ h\in H_{\mathbb R}: sign
(\chi_{\alpha_i}(h))=\epsilon_i \text { for all } \alpha_i\in \Pi \}$, and we have:
$$H_{\mathbb R}=\underset {\epsilon\in {\mathcal E}}  {\bigcup}H_{\epsilon} .$$

\begin{notation}
\label{D2.2.1}
We need notation to parametrize connected the components of a split Cartan subgroup, roots and root characters. Unfortunately we need such notation for {\em all} the parabolic subgroups associated to arbitrary  subsets of $\Pi$. Recall (Notation  \ref{D2.1.1}) that notation associated to a parabolic subgroup determined by a subset $\Pi\setminus A$, $A\subset \Pi$ is usually indicated by changing the standard notation with the use of a superscript ${\cdots }^A$. Thus we have : $\Delta^A\subset \Delta$, root system giving rise to a {\em semisimple} Lie algebra ${\mathfrak l}^A \subset {\mathfrak g}$. Also there are corresponding connected  semisimple Lie subgroups $ L_{\mathbb C}^A \subset G_{\mathbb C}$, $ L^A \subset  L_{\mathbb C}^A$. The adjoint group is denoted by  $L_{\mathbb C}(\Delta^A)=Ad(L_{\mathbb C}^A)$ and  it has a real connected Lie subgroup $L(\Delta^A)$.  Let ${\tilde L}(\Delta^A)$ be defined in the same way as $\tilde G$ but relative to the  root system $\Delta^{A}$.  We let ${\mathfrak h}^A$ be the real span of the $h_{\alpha_i} $ with $\alpha_i
\not\in A$. This is a (split) Cartan subalgebra of ${\mathfrak
l}^A$ denoted as $H_{\mathbb R}^A$.  The corresponding connected Lie subgroup is $H^A=\exp({\mathfrak h}^A)$ (exponentiation taking place inside $G_{\mathbb
C}$). 
\end{notation}

We also consider Lie subgroups of $\tilde G$ corresponding to the subsystems
of the Toda lattice.  In accordance to our convention for Levi subgroups associated with $A\in \Pi$, we have  
$${\mathcal E}^A=\{\ \epsilon
=(\epsilon_1,..,\epsilon_l): \epsilon_i=1 \text { if } \alpha_i \in A \}.$$
 Thus $H^A$ is by definition a subgroup of $H$. This  Lie subgroups of $H$  is
isomorphic to the connected component of the identity of a Cartan subgroup of a
real semisimple   Lie group  that corresponds to ${\mathfrak l}^A$.
There is a bijection 
\begin{equation}
{\mathcal E}^A \cong { \mathcal E}(\Delta^A)=\left\{
(\epsilon_{j_1},..,\epsilon_{j_m}) :\ 
\alpha_{j_i}\not\in A  \text { for } i=1,\ldots,m=|\Pi\setminus A| \right\}
 \label{2.2.3}
\end{equation}
which is given in the obvious way by restricting a function
$\epsilon :
\Pi \to \{\ \pm 1 \}$ such that $\epsilon_i=\epsilon(\alpha_i)=1$ for
$\alpha_i\in A$ to a new function $\epsilon (\Delta^A)$ with domain  $\Pi\setminus A$.

We now consider the Cartan subalgebras and the Cartan subgroups for Levi factors of parabolic subalgebras and subgroups dtermined by $A\subset \Pi$:
\begin{defn}
\label{D2.2.2}
For $A\subset \Pi$, we denote:
\begin{itemize}
\item  $H_{\mathbb R}^A =\underset{\epsilon \in {\mathcal E}^A }
{\bigcup }h_\epsilon H^A$ (if $\Delta^A=\emptyset$
(i.e. $A=\Pi$) then $H_{\mathbb R}^A=\{e\}$).

\item  $H_{\mathbb R} ( \Delta^A)$ a Cartan
subgroup of  ${\tilde L}(\Delta^A)$, defined in the same way as $H_{\mathbb R}$ and
having Lie algebra ${\mathfrak h}^A$ ($H_{\mathbb R}(\Delta^A)=\{e \}$ if
$\Delta^A=\emptyset$).

\item   $H_\epsilon^{A,\leq}=\{\ h\in
h_\epsilon H^A:\text { for all } \alpha_i \in \Pi\setminus A: |\chi_{\alpha_i}(h)| \leq 1
\}$  the antidominant chamber of $H^A_\epsilon=h_{\epsilon}H^A$. Similarly we
consider $H_\epsilon^{A, <}$ using {\it strict} inequalities.

\item  $H(\Delta^A)_\epsilon^{\leq }=\{\ h\in H(\Delta^A)_\epsilon:
\text { for all } \alpha_i \in \Pi\setminus A: |\chi^{\Delta^A}_{\alpha_i}(h)| \leq 1 \}$.
Similarly we consider a version with strict inequalities.

\item $\chi^{\Delta^A}_{\alpha_i}$  the root character associated to $\alpha_i$
on the Cartan $H_{\mathbb R}(\Delta^A)$.
 \end{itemize}
We use notation $\cdots^{\leq }$ to indicate the antidominant chamber on Cartan subalgebras and subgroups and  $\cdots^{< }$ for strictly antidominant chambers.
\end{defn}

\begin{example}
\label{E2.2.2}
In the case of $\tilde G=SL(3, {\mathbb R})$, 
$H_{\mathbb R}$ is the group
\begin{equation}
H_{\mathbb R}=\left\{\ {\rm diag} \left(a, b,c\right)\ : \ a\not=0, ~b\not=0, ~abc=1~\right\}.
\nonumber 
\end{equation}
For $A=\{ \alpha_1 \}$, $H^A$ is the group $\{\ {\rm diag} (1,
a,a^{-1}): a> 0 \}$. The set ${\mathcal E}^A$ consists of $(1,1)$ and $(1,-1)$.
The element $ h_{(1,-1)}={\rm diag} (-1,-1,1)$ and thus $h_{(1,-1)} H^A=\{\ {\rm diag} (-1, -a,a^{-1}) : a> 0 \}$. These two components form $H^A_{\mathbb R}$. 
Note that $H_{\mathbb R} L^A$ is  the Lie subgroup of
$SL(3, {\mathbb R})$ consisting of all real matrices of the form,
\begin{equation}
H_{\mathbb R}L^{\{\alpha_1\}}=
\pmatrix a& 0& 0
\cr 0& b& c   \cr 0& d& e  \cr \endpmatrix
\nonumber
\end{equation}
having determinant one.
The group $L^A$ is obtained by setting $a=1$ and the determinant equal to one.
The group $L(\Delta^A)$ is isomorphic to $Ad(SL(2,{\mathbb R}))$  the adjoint group
obtained from $L^A$ and $\tilde L(\Delta^A)$ is isomorphic
$Ad(SL(2,{\mathbb R})^{\pm})$. Thus these three groups $L^A$, $L(\Delta^A)$
and $\tilde L(\Delta^A)$ are all different in this case.
\end{example}

\begin{defn}
\label{D2.2.3}
Recall that the $S$ is a  subset of $\Pi$ indicating {\em colored} vertices in a Dynkin diagram. Let $\eta: S\to \{ \pm 1 \}$ be any function. We let
 $\epsilon_\eta\in {\mathcal E}^{\Pi\setminus S}$ defined by
$\epsilon_\eta (\alpha_i)= \eta(\alpha_i)$ if $\alpha_i\in S$, $\epsilon_\eta
(\alpha_i)=1$ if $\alpha_i\not\in S.$
Thus there is a bijective correspondence between the set of all  functions
$\eta$
and the set ${\mathcal E}^{\Pi\setminus S}$ and (by equation (\ref{2.2.3})) a second bijection
with the set
${\mathcal E}(\Delta^{\Pi\setminus S} )$: 
\begin{equation}
\eta\mapsto \epsilon_\eta \in {\mathcal E}^{\Pi\setminus S},
\nonumber
\end{equation}
and  
\begin{equation}
\eta\mapsto \epsilon_\eta (\Delta^{\Pi\setminus S}) \in {\mathcal E}
(\Delta^{\Pi\setminus S}).
\nonumber
\end{equation}
The  exponential map $h_\epsilon \exp :\mathfrak h \to h_\epsilon H$ allows us to define chamber walls in $H$ and  therefore in  any of the connected components of $H_{\mathbb R}$. For any root $\phi \in \Delta$ the set $\{\ h\in h_\epsilon H:
|\chi_{\phi}(h)|=1 \}$ defines the  $\phi$-wall of $H_\epsilon=
h_\epsilon H$.
The intersection of all the $\alpha_i$-walls of $H_\epsilon$ is the set $\{ h_\epsilon \}$. This is also the intersection of all the
$\phi$ walls of $H_\epsilon$ with  $\phi \in \Delta$.
\end{defn}

We also define the following subsets of $H_{\epsilon}$:
\begin{defn}
\label{D2.2.4}
We denote:
\begin{itemize}
 \item  ${\mathcal D}={\mathcal D}(\Delta)=\{\ h_\epsilon:~ \epsilon \in  {\mathcal E}~ \}$.
 \item  ${\mathcal D} (\Delta^{\Pi\setminus S})=\{\ h_\epsilon: ~\epsilon \in {\mathcal E}(\Delta^{\Pi\setminus S})~\}$ .
\end{itemize}
\end{defn}

The set ${\mathcal D}$ has two structures: it is a finite group and also a set
with an
action of $W$. In Proposition \ref{P2.2.2} it is the first structure that is
emphasized but in the proof of Proposition  \ref{P3.3.1}
it is the second structure which is
relevant. We now look at the $W$ action.

Since $W$ acts on $H_{\mathbb R}$ and $w\in W$ sends a $\phi$-wall of
$H_\epsilon$ to
the $w(\phi)$-wall of some other $H_{\epsilon^\prime }$, this set  ${\mathcal D}$ is
preserved by $W$ and thus acquires a $W$ action. Given $S\subset \Pi$   we
similarly obtain that ${\mathcal D}(\Delta^{\Pi\setminus S})$ has a $W_S$ action.

The map 
\begin{equation}
{\mathcal D} (\Delta^{\Pi\setminus S}) \to {\mathcal E} (\Delta^{\Pi\setminus S})
 \label{2.2.8}
\end{equation}
sending $h_\epsilon \mapsto \epsilon$ also defines  an action of $W_S$
on the set of signs ${\mathcal E}(\Delta^{\Pi\setminus S}).$  Recall (Notation  \ref{D2.2.1}) that
${\mathcal E}^{\Pi\setminus S} \subset {\mathcal E}$ denotes those $\epsilon$ for which
$\epsilon_i=1$ whenever  $\alpha_i\not\in S$. We have   a bijection (by (\ref{2.2.8}) together with
(\ref{2.2.3}) )
 \begin{equation}
{\mathcal E}^{\Pi\setminus S} \cong {\mathcal D}(\Delta^{\Pi\setminus S}).
\nonumber
\end{equation}
 The
$W_S$
action on the set ${\mathcal D}(\Delta^{\Pi\setminus S})$  thus gives a $W_S$ action on ${\mathcal
E}^{\Pi\setminus S}$ such that for any $\epsilon\in {\mathcal E}^{\Pi\setminus S}$ only the $\epsilon_j$
with $\alpha_j\in S$ may change in sign  under the action. Note that this
construction requires looking at the $h_\epsilon$ with $\epsilon \in {\mathcal
E}(\Delta^{\Pi \setminus S})$ in the adjoint representation of ${\mathfrak l}^A$.

\begin{rem}
\label{R2.2.1}
The root characters can be expressed as a product
of the
simple root characters raised to certain integral powers $\phi= { \sum_{i=1}^l } c_i\alpha_i \,$ with $c_i\in \mathbb Z$, $e^{\phi }=
 { \prod_{i=1}^l } \chi_{\alpha_i}^{c_i} \,  $.
Therefore if
$h\in H_{\mathbb C}$ the scalars $\chi_{\alpha_i}(h)$ determine all the scalars
$e^\phi (h)=\chi_\phi (h)$ and thus $h$ is uniquely determined. Moreover
$e^{\langle \phi, x+\sqrt {-1} y\rangle }$ is real for all $\phi \in \Delta$ if and
only if
$e^{\langle \alpha_i, x+\sqrt {-1} y\rangle }$ is real for all $\alpha_i \in \Pi$.
\end{rem}

\begin{prop}
The Cartan subgroup $H_{\mathbb R}$ of ${\tilde G}$ has
$2^l$ components. We have $H_{\mathbb R} = {\mathcal D} H$ where ${\mathcal D}$ (Definition
\ref{D2.2.4})  is the finite group of all the $h_\epsilon$, $\epsilon\in {\mathcal E}$.\label{P2.2.2}
\end{prop}
\begin{Proof}
 It was shown in Lemma  \ref{L2.2.1} that  the elements
$\exp(x+\sqrt {-1} y)$
with $x,y\in \mathfrak h$ such that $e^{\langle \alpha_i ; x+\sqrt {-1} y\rangle }$ is
real for all
$\alpha_i
\in \Pi$ are those for which $y$ has the form $y=\sum_{i=1}^{l}k_i y_i$ with
$y_i$ as in Notation \ref{D2.1.2}, and $k_i$ are integers. Moreover as in
Remark  \ref{R2.2.1}  we
also
have that $e^{\langle \phi , x+\sqrt {-1} y\rangle } $ is real for any $\phi\in
\Delta$ exactly when
$y=\sum_{i=1}^{l} k_iy_i$ for some integers $k_i$ . Therefore all the root
spaces of $\mathfrak g$ are stabilized under the adjoint action of $\exp(x+\sqrt {-1} y)$.
Since
clearly $\mathfrak h$ is also stabilized, then $\exp(x+\sqrt {-1} y)$ stabilizes all of
$\mathfrak
g$ and this implies that $\exp(x+\sqrt {-1} y)\in {\tilde G}$. In fact this shows that
$\exp(x+\sqrt {-1} y) \in \tilde G$ if and only if $y$ has the form
$y=\sum_{i=1}^{l} k_i
y_i$ for certain $k_i$ integers. Thus Lemma \ref{L2.2.1} and these remarks
compute the
intersection $\tilde G\cap H_{\mathbb C}$. From here it is easy to conclude.
\end{Proof}

Note that Proposition \ref{P2.2.2} implies that ${\mathcal D} (\Delta) $ in Definition
\ref{D2.2.4}  is isomorphic to  $H_{\mathbb R}/ H$ as a set with a $W$ action. That
the $W$ actions agree is verified in Corollary \ref{C3.3.1}.

\vskip 0.3cm

We now give the proof of Proposition \ref{P2.2.1}:

\begin{Proof}
Let $Ad$ denote the representation of
$GL(n,\mathbb C)$ on ${\mathfrak{sl}}(n,{\mathbb C})$.
Then we have $Ad(GL(n, {\mathbb C}))=Ad(SL(n, {\mathbb C}))$. If $n$ is odd, $Ad(SL(n, {\mathbb C} ))$ is isomorphic to $SL(n,{\mathbb C})$.
Denote
$D_i=diag (r_1,..,r_n)$ with $r_i$ as in Example \ref{E2.2.1}. we set $\bar h_i=D_i$
when $n$ is even and $\bar h_i=(-1)^{n-i} D_i$ when $n$ is odd. Let $h_i=
Ad(\bar h_i)$.  We have that for  each $h_i\in {\tilde G}$,
$\chi_{\alpha_i}(h_i)=e^{\pi \sqrt {-1} \delta_{i,j}}$ for $\alpha_i\in \Pi$.
The $h_i$ generate the group ${\mathcal D}$ with $2^l$
elements where $l=n-1$ and now the group $G_1=\langle Ad(SL(n, {\mathbb R})),
h_i, i=1,..,l\rangle$ is a subgroup of $\tilde G$ and it is isomorphic to
$ Ad(SL(n, {\mathbb R})^{\pm})$ for
$n$ even and to $Ad(SL(n, {\mathbb R}))\cong SL(n, {\mathbb R})$ for $n $ odd.  What
remains is to verify that  ${\tilde G} \subset G_1$.

Let ${\tilde K}=\{ Ad(g): g\in U(n) \}\cap { \tilde G}$. By the Iwasawa
decomposition of ${\tilde G}$ and $G_1$ , it suffices to show that
${\tilde K}=K {\mathcal D}$, $K=SO(n)$. The right side, $K {\mathcal D}$ is either $O(n)$ or
$SO(n)$ according to the parity of $n$.
Recall that the maximal compact Lie subgroup
 $\tilde K$ of $\tilde G$ acts transitively on the
set $X$ of all maximal abelian Lie subalgebras $\mathfrak a$ which are
contained in the vector space $\mathfrak p$. The action of $K=Ad(SO(n))$
is also transitive on $X$ (by (2.1.9) of \cite{WA}) and thus for any $g\in
\tilde K$ there is $k\in K$
such that $g=kD$ where $D$ is the isotropy group (in $\tilde K$) of an
element in $X$,
for instance of the element ${\mathfrak h} \in X$. However
this isotropy group $D$ has been computed implicitly in the proof of
Proposition \ref{P2.2.2} and $D={\mathcal D}$. Therefore ${\tilde K}=K{\mathcal D}$.
\end{Proof}

\begin{prop}
Let $\alpha_i\in \Pi$ and assume that
$\epsilon=(\epsilon_1,..,\epsilon_l) \in {\mathcal E}$. Then
$s_{\alpha_i}h_\epsilon=h_{\epsilon^\prime} $ where
$\epsilon^\prime_j=\epsilon_j
\epsilon_i ^{-C_{j,i}}$ with $(C_{i,j})$ the Cartan matrix. \label{P2.2.3}
\end{prop}
\begin{Proof}
 This  follows from the expression of the Weyl group action on
elements in ${\mathfrak h}_{\mathbb C}^\prime $ given by: $s_{\alpha_i}x=x-(\check
\alpha_i,
x) \alpha_i$. If  this expression is applied to $x=\alpha_j$ it gives
$s_{\alpha_i} \alpha_j=\alpha_j- C_{j,i}\alpha_i$. On the level of root
characters this becomes, by exponentiation of the previous identity,
\begin{equation}
s_{\alpha_i} \chi_{\alpha_j}= \chi_{\alpha_j}\chi_{\alpha_i}^{-C_{j,i}}.
\nonumber
\end{equation}

 Now recall that $\epsilon_j$ is just $\chi_{\alpha_j}$ evaluated at
$h_\epsilon$. Also $\epsilon^\prime_j$ will be $\chi_{\alpha_j}$ evaluated at
$s_{\alpha_i}h_\epsilon$. When we evaluate $\chi_{\alpha_j}$ on
$s_{\alpha_i}h_\epsilon $ in order to compute the corresponding j-th sign, we
obtain  $\chi_{s_{\alpha_i} \alpha_j }(h_\epsilon)$. Therefore  the sign of
$\chi_{\alpha_j}$ on the $s_{\alpha_i}h_\epsilon$ is given by the product
$\epsilon_j \epsilon_i^{-C_{j,i}}$. Finally we use the fact that the set of all
scalars $\chi_{\alpha_i}(h)$ determines $h$. Thus $\epsilon^\prime$ determines
the element $h_{\epsilon^\prime}$ giving rise to the equation
$s_{\alpha_i}h_\epsilon= h_{\epsilon^\prime}$.
\end{Proof}

The sign change $\epsilon_j \to \epsilon_j'$ in Proposition \ref{P2.2.3}
is precisely the gluing rule in Lemma 4.2 for the indefinite Toda
lattice in \cite{KO-YE2}.  Then the gluing pattern using the Toda
dynamics is just to identify each piece of the connected component
$H_{\epsilon}$ (see Figure \ref{F1.1.1}). The sign change on subsystem corresponding to $H_{\mathbb R}^A$ with $A\subset \Pi$ can be also formulated as:

\begin{prop}
Let $\alpha_i\in \Pi\setminus A$ and assume that
$\epsilon=(\epsilon_1,..,\epsilon_l) \in {\mathcal E}^A$. Then

a)  if $\epsilon_i=1$, $s_{\alpha_i}h_\epsilon =  h_\epsilon$. If
$\epsilon_i=-1$
 then $s_{\alpha_i} h_\epsilon = h_{\epsilon^\prime}$ where
$\epsilon^\prime_j=\epsilon_j \epsilon_i^{-C_{j,i}}$. In addition
$h_{\epsilon^\prime}$ factors  as a product $( \prod_{\alpha_j\in A,C_{j,i}
\text
{ is odd }  } h_j ) h_{\epsilon_A}$ where $\epsilon_A\in {\mathcal E}^A$.

b)  if $\epsilon_i=1$, $s_{\alpha_i} H_\epsilon^A \subset H_\epsilon^A$. If
$\epsilon_i=-1$  then $s_{\alpha_i} H_\epsilon^A \subset
h_{\epsilon^\prime}H^A$
where $\epsilon^\prime_j=\epsilon_j \epsilon_i^{-C_{j,i}}$. In addition
$h_{\epsilon^\prime}$ factors  as a product $( \prod_{\alpha_j\in A,C_{j,i}
\text
{ is odd }  } h_j ) h_{\epsilon_A}$ where $\epsilon_A\in {\mathcal E}^A$.

c) The sign $\epsilon_k= sign ( \chi_{\alpha_k}(h) )$ for  any $h\in
H_\epsilon$
agrees with  $sign (\chi_{s_{\alpha_i}\alpha_k }(h^\prime))$ for any $h^\prime
\in s_{\alpha_i}( H_\epsilon )$.
 \label{P2.2.4}
\end{prop}

\begin{Proof} Part a) follows from Proposition \ref{P2.2.3} but with the observation
that in, this case,  $\epsilon^\prime $ may fail to be in ${\mathcal E}^A$ even if
$\epsilon\in {\mathcal E}^A$. This happens exactly  when  $\epsilon_i=-1$ and
$\alpha_j\in A$ with $C_{j,i}$ odd (either $-3$ or $-1$).Under these
circumstances  $\epsilon_j^\prime=-1$ (rather than one as required in the
definition of ${\mathcal E}^A$). We can fix this problem by factoring
$h_{\epsilon^\prime}$ as a product $( \prod_{\alpha_j\in A,C_{j,i} \text {
is odd
}  } h_j)h_{\epsilon_A}$ where $\epsilon_A\in {\mathcal E}^A$.

Part b) follows from part a) and the fact that each $h_\epsilon H$ is a
connected
component of $H_{\mathbb R}$ in Proposition \ref{P2.2.2}.

Part c) follows easily from
$\chi_{s_{\alpha_i}\alpha_k}(s_{\alpha_i}h)=\chi_{\alpha_k}(s_{\alpha_i}
(s_{\alpha_i}h) )= \chi_{\alpha_k}(h)$. Also by the fact that
 the sign of a root
character is
constant along a connected component and that we have $h\in H_\epsilon $ and
$s_{\alpha_i}h$ is in the new component $s_{\alpha_i} H_\epsilon$. The two
desired signs have thus been computed in the two connected components and they
agree.
\end{Proof}
\begin{rem}
From  Proposition \ref{P2.2.4} it follows that any
connected component $H_\epsilon^A,$ is the union of chambers of the form
$w\left(H_{\epsilon (w)}^{A,\leq}\right)$,with  $w\in W_{\Pi\setminus A}$, $\epsilon (w) \in {\mathcal E}$,
\begin{equation}
H^A_{\epsilon}=\bigcup_{w\in W_{\Pi\setminus A}}w\left(H^{A,\leq }_{\epsilon (w)}\right).
\nonumber
\end{equation}
\label{R2.2.2}
\end{rem}

\section{Colored Dynkin diagrams}
 \label{S3}

We now introduce some notation that will ultimately  parametrize the cells in a
cellular decomposition of the  smooth compact manifold $\hat H_{\mathbb R}$ to be defined in \S \ref{S5}.

\subsection{Colored Dynkin Diagrams}
 \label{ss3.1} 

Let us first define:
\begin{defn} \label{D3.1.1}
 {\it Colored Dynkin diagrams $\mathbb D(S)$: } 
A colored Dynkin diagram is a Dynkin diagram where all the
vertices
in a set $S\subset \Pi$ have been colored either red $R$ or blue $B$.
For example, in ${\mathfrak{sl}}(4, {\mathbb R})$,  $\circ_R-\circ -\circ_B $ is a
colored Dynkin
diagram with $S=\{ \alpha_1, \ \alpha_3 \}$.
 Thus a colored Dynkin diagram where $S\not=\emptyset$, corresponds to a pair
$(S,\epsilon_\eta)$ with $S\subset \Pi$ and  $\eta :S\to \{ \pm 1 \}$ any
function. Here $\eta (\alpha_i)=-1$ if $\alpha_i$ is colored $R$ and $\eta
(\alpha_i)=1$ if $\alpha_i$ is colored $B$. If $S=\emptyset$ then $\epsilon_o$
with $\epsilon_o (\alpha_i)=1$ for all $\alpha_i\in \Pi$ replaces
$\epsilon_\eta$. We denote:
\begin{itemize}
\item $D=(S,\epsilon_{\eta})$ or $(S, \epsilon_{\eta}(\Delta^{\Pi\setminus S}))$
with $\epsilon_{\eta} \in {\mathcal E}^{\Pi\setminus S}$
\item  ${\mathbb D}(S)=\{\  D=(S,\epsilon_{\eta})~ :\text {
the vertices in} S \text { are colored}\} $
\end{itemize}
We also introduce {\it an oriented colored Dynkin diagram} which is defined as  
 a pair $(D, o)$ with $o \in \{\pm 1 \}$ and $D$ a
colored
Dynkin diagram.
\end{defn}

\subsection
{Boundary of a colored Dynkin diagram} \label{ss3.2}
We now define the boundaries of a cell parametrized by a colored Dynkin diagram $D$.
\begin{defn} {\it The boundary $\partial_{j,c} D$:}
For each $(j,c)$ with $c=1,2$ and $j=1,..,m$ we define  a new colored Dynkin
diagram $\partial_{j,c}D$, the $(j,c)$ -boundary of the  $D$  by
considering $\{
\alpha_{i_j}: 1\leq i_1 < ...< i_m\leq l \}$ the set $\Pi\setminus S$ of uncolored
vertices and $m=|\Pi\setminus S|$. The
$\partial_{j,c} D$ is then a  new colored Dynkin diagram obtained by
coloring the
$i_j$ th vertex with $R$ if $c=1$ and with $B$ if $c=2$. The boundary  of an
oriented colored Dynkin diagram $(D,o)$, $o\in \{\ \pm 1 \}$ is, in
addition, given an orientation defined to be  the sign $(-1)^{j+c+1}o$. Recall that a
colored Dynkin diagram $D$ corresponds to a pair $(S,\epsilon_\eta)$ with
$S\subset \Pi$ and $\eta:S\to \{\ \pm 1 \}$ (Definition \ref{D2.2.3}).
Thus the boundary $\partial_{j,c}$ determines a new pair $(S\cup \{\
\alpha_{i_j}
\}, \epsilon_{\eta^\prime} )$ associated to $\partial_{j,c}D$.
We can then define the following  boundary maps,
\begin{equation}
(-1)^{j+c+1}  \partial_{j,c}:  {\mathbb Z}\left[{\mathbb D}(S)\right]\, \to
 {\mathbb Z}\left[{\mathbb D}(S\cup \{\ \alpha_{i_j} \})\right] \, 
\label{3.2.4}
\end{equation}
\end{defn}

\begin{figure}
\includegraphics{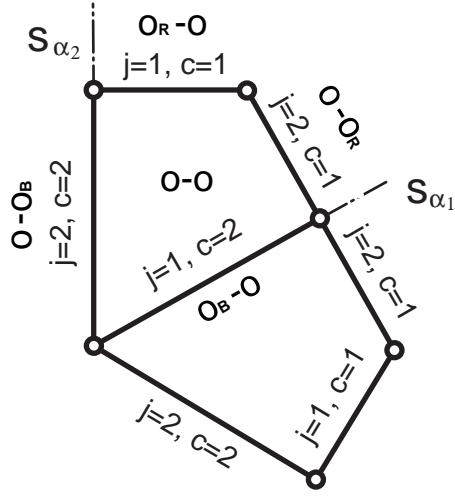}
\caption{The boundary $\partial_{j,c}$ in the case of  ${\mathfrak {sl}}(3, \mathbb R)$.}
\label{Fboundary}
\end{figure}

\begin{example}
\label{E3.2.1} 
The boundary of  $\circ  - \circ $ (which we can
picture as a box) consists of segments (one dimensional boxes) given by
$\partial_{1,1}(\circ-\circ)=\circ_R-\circ $,
$\partial_{2,1}(\circ-\circ)=\circ  - \circ_R$,
$\partial_{1,2}(\circ-\circ)=\circ _B -\circ $ and
$\partial_{2,2}(\circ-\circ)=\circ -\circ _B$.
The orientation sign associated to $\circ_R-\circ $ is the following: with
$c=1$ and $j=1$ one has  $(-1)^{j+c+1}=(-1)^{3}=-1$. 
We illustrate the example in Figure\ref{Fboundary}.
\end{example}

\subsection {$W_S$-action on  colored Dynkin diagrams}
 \label{ss3.3} 

We now move these colored Dynkin diagrams around with elements in $W$.  
A $W_S$-action on the diagram $D\in{\mathbb D}(S)$, $W_S:{\mathbb D}(S)\to{\mathbb D}(S)$, is defined as follows:
\begin{defn}
\label{D3.3.1}
 For any $\alpha_i\in S$ (the $\alpha_i$ vertex is colored),
$s_{\alpha_i}D=D'$ is a new colored Dynkin diagram having the colors according
to the sign change $\epsilon_j'=\epsilon_j\epsilon_i^{-C_{j,i}}$ in
Proposition \ref{P2.2.4}
with the identification that $R$ if the sign is $-1$, and $B$ if it is $+1$.
For example, in the case of ${\mathfrak{sl}}(3, {\mathbb R})$,
$s_{\alpha_1}(\circ_R-\circ_B)=\circ_R-\circ_R,
s_{\alpha_1}(\circ_B-\circ_R)=\circ_B-\circ_R$.

We also define a $W_S$ action on the set ${\mathbb D} (S)\times \{\pm 1 \}$ of
oriented colored Dynkin diagrams. If $\alpha_i\in S$, the action of
$s_{\alpha_i}$ on the pair $(D,o)$ is given by
$s_{\alpha_i}(D,o)=(s_{\alpha_i}D,
(\epsilon_i)^{r_{\alpha_i} }o)$ where $ r_{\alpha_i}$ is the number of elements
in the
set $\{\ \alpha_j\in \Pi\setminus S : C_{j,i} \text { is odd } \}$ and $\epsilon_i=\pm 1$ depending on the color  of $\alpha_i$.
\end{defn}
We confirm the Definition:
\begin{prop}
Definition \ref{D3.3.1} above gives a well-defined
action of $W_S$ on the set ${\mathbb D}(S)$ of  colored Dynkin diagrams with set of
colored vertices $S$. \label{P3.3.1}
\end{prop}

 \begin{Proof} This follows from Proposition \ref{P2.2.3} and the correspondence
 $(S,\epsilon_\eta) \to h_{\epsilon_\eta}$ giving a bijection between
 ${\mathbb D}(S)$ and ${\mathcal D}(\Delta^{\Pi \setminus S})$.
\end{Proof}

\begin{cor}
There is a bijection of  ${\mathcal D}(\Delta ) \to
H_{\mathbb R}/ H$ intertwining the $W$ actions on both sets. \label{C3.3.1}
\end{cor}

\begin{Proof}
The two sets ${\mathcal D}(\Delta)$ and $H_{\mathbb R}/H$ are clearly in
bijective correspondence and an element $h_\epsilon$ corresponds to the coset
$h_\epsilon H$.  By Proposition  \ref{P2.2.4},  Notation \ref{D3.3.1} and
Proposition \ref{P3.3.1}, the actions on these two sets agree.  They are both
given by Part $a$) in
Proposition \ref{P2.2.4}.
\end{Proof}

\begin{notation}
\label{N3.3.1}
Motivated by Example  \ref{E5.2.2} we consider the set $ W \underset { W_S }
   {\times  } {\mathbb D}(S) \,$ the $W$ translations of the set  ${\mathbb D}(S)$. We
write the
elements in this  set as pairs $(w,D)$ and introduce an equivalence relation
$\sim$ on these pairs,   where for any $x\in W_S$, $(wx,D)\sim (w,xD)$. Denote
the equivalence classes $[(w,D)]$. This set of equivalence classes  is then  in
bijective correspondence with the set   ${\mathbb D}(S)\times W/W_S$ of pairs $(D,
[w])$ with $[w]=[w]^{\Pi\setminus S}\in W/W_S$ and $D$ in ${\mathbb D}(S)$. The
correspondence
is such that $[(w^\bullet,D)]$ corresponds to $(D, [w^\bullet])$ with
$w^\bullet
\in [w]^{\Pi\setminus S}$ a minimal length representative of a coset in $W/W_S$.

We denote
\begin{equation}
{\mathbb D}^k=\{\ (D, [w]^{\Pi\setminus S}) : S\subset \Pi, |S|=k, w\in W \},
\nonumber
\end{equation}
which also parametrizes all the connected components of the Cartan
subgroups of the form $H_{\mathbb R}(w(\Delta^{\Pi\setminus S}))$.
\end{notation}

\begin{rem} For the reader who is not interested in the $W$  action or torsion,
 the object defined in  Definition \ref{D3.3.2} becomes over $\mathbb Q$ the vector space with 
 basis given by the full set of colored Dynkin diagrams having a fixed number of uncolored
  vertices.  The boundary maps are obtained by {\em translating} around the $\partial_{j,c}$
   defined for the antidominant chamber (see (\ref{3.3.9}) below). Perhaps 
   this boundary construction can be appreciated in Example \ref {E5.2.2} in subsection 
   \ref{SSE5.2.2} where portions of the boundary of a chamber are translated and some identifications take place. The tensor product notation below accomplishes the required boundary identifications algebraically.
\label{Rrational}
\end{rem}

\begin{defn}\label{D3.3.2}
 {\it The ${\mathbb
Z}[W]$-modules ${\mathcal M}(S)$:} 
The full set of colored Dynkin diagrams is the set  $\mathbb
D(S)\times W/W_S$ of all pairs $(D, [w]^{\Pi\setminus S})$.

We can also define the full set of  $oriented$ colored Dynkin diagrams by
considering ${\mathbb D}(S)\times \{ \pm 1 \} \times W/W_S$ for different subsets
$S\subset \Pi$. As $W$ sets these correspond to  $W \underset { W_S }  {
\times } {\mathbb D} (S)\times \{ \pm 1 \} \,$.
\end{defn}

If we consider ${\mathbb D}(S)\times \{ \pm 1 \}$ as imbedded in ${\mathbb Z} [{ \mathbb
D}(S)]$ by sending $(D, o)$ to $oD$, $o \in \{\pm 1\}$ we may consider a $W_S$-action on $\pm { \mathbb D}(S)$, namely the action on oriented colored Dynkin
diagrams.
This produces a ${\mathbb Z} [W]$-module ${\mathbb Z}[ W] \underset { {\mathbb Z} [W_S ] }
{ \otimes } {\mathbb Z} [{\mathbb D}(S)] \,$, and we denote this module by
\begin{equation}
{ \mathcal M}(S)={\mathbb Z} [W]
\underset { {\mathbb Z}[W_S] } { \otimes }{ \mathbb Z} [{\mathbb D}(S)] \,.
\nonumber
\end{equation}
Also denote by ${\mathcal M}_{l-k}$ the direct sum of all these modules over all
sets $S$
with exactly $k$ elements,
\begin{equation}
{\mathcal M}_{l-k}=\bigoplus_{|S|=k}{\mathcal M}(S)\ .
\nonumber
\end{equation}

\begin{rem}
\label{R3.3.1}
Assume that  $A$ is a ${\mathbb Z}[W_S]$-module and $B$
is  a
${\mathbb Z}[W_{S^\prime}]$-module with $S\subset S^\prime$. In particular $B$
can be
regarded as a ${\mathbb Z}[W_S]$-module by restriction.  Let $f:A\to B$ be a map
intertwining these two  ${\mathbb Z}[W_S]$-module structures involved. Then,
tensoring
with ${\mathbb Z}[W]$ we obtain a map of ${\mathbb Z}[W]$-modules:
\begin{equation}
F(f):{\mathbb Z}[W] \underset { \mathbb Z[W_S] }   {
\otimes } A \, \to \mathbb Z[W] \underset { \mathbb Z[W_S] }   { \otimes } B \ .
\nonumber 
\end{equation}
Since, in addition, $B$ is a
${\mathbb Z}[W_{S^\prime }]$-module where $S\subset S^\prime$
then there is a second map:
\begin{equation}
g: {\mathbb Z}[W] \underset { {\mathbb Z}[W_S] }   { \otimes } B \ .  \to{\mathbb Z}[W]
\underset {{ \mathbb Z}[W_{S^\prime}] }     { \otimes } B \,
\nonumber
\end{equation}
Let $g\circ F(f)=T(f)$. This is then a map of ${\mathbb Z}[W]$ modules
\begin{equation}
T(f): {\mathbb Z}[W] \underset { {\mathbb Z}[W_S] } { \otimes } A \,\to  {\mathbb
Z}[W] \underset { {\mathbb Z}[W_{S^\prime}] }   { \otimes } B \ . 
\nonumber
\end{equation}
\end{rem}

We now give the boundary maps of ${\mathcal M}_{l-k}$.
 We regard ${\mathbb D}(S)\times \{ \pm  1\}$ as a
${\mathbb Z}[W_S]$-module using
Definition \ref{D3.3.1} the oriented case.  We can view a pair $(D, o)$
instead as $\pm D\in {\mathbb Z} [{\mathbb D}(S)]$.
This gives a ${\mathbb Z}[W_S]$-module structure to each of
the $\mathbb Z$ modules involved on the domain of the map in  (\ref{3.2.4})
 and a $W_{S\cup \{\ \alpha_{i_j} \} }$ to those on the
co-domain. The map in (\ref{3.2.4})  now has the form
described in Remark \ref{R3.3.1}. We apply the construction $T(\partial_{j,c} )$
in Remark \ref{R3.3.1}  to (\ref{3.2.4})  and add over all possible subsets $S\subset \Pi$ with $|S|=k$ on the left side and with $|S|=k+1$ on the right side.
We then obtain the boundary maps,
\begin{equation}
\partial_{l-k}: {\mathcal M}_{l-k}\to {\mathcal M}_{l-(k+1)}~,
 \label{3.3.6}
\end{equation}
which are all given by 
\begin{equation}
\partial_{l-k}( w^\bullet
\otimes X)
=\displaystyle{\sum_{j=1}^{l-k}\sum_{c=1}^2 (-1)^{j+c+1} w^\bullet
 \otimes T(\partial_{j,c}) X} ~,
\label{3.3.9}
\end{equation}
where $X\in {\mathbb D}(S)$ and
$w^\bullet \in [w]^{\Pi\setminus S}$ is the minimum length representative in $W/W_S$ (see Notation \ref{N3.3.1}).
We thus have:

\begin{prop}
The maps $\partial_{l-k}$ of (\ref{3.3.9}) define
a chain
complex ${\mathcal M}_*$ of ${\mathbb Z}[W]$ modules. \label{P3.3.2}
\end{prop}

 \section{Cartan subgroups  and  Weyl group actions}
  \label{S4}

We here discuss relations between some Cartan subgroups of Levi factors, and
verify the $W_S$-action on oriented colored Dynkin diagrams.

\subsection{ Cartan subgroups of Levi factors} \label{ss4.1} 

 Let $Ad^{\Delta^A}$ denote the adjoint
representation of the Lie subgroup $ H_{\mathbb R} L^A$ of $G$ on the  Lie algebra
${\mathfrak l}^A$ . Note that we are deviating slightly from the standard convention  and
denoting  by
$Ad^{\Delta^A}$
the  representation of $H_{\mathbb R} L^A$  acting on the semisimple part of the Levi factor ${\mathfrak
l}^A$, seen as a {\it quotient} of  the group  action  on  ${\mathfrak h}^A + {\mathfrak
l}^A$ (thus dividing by the center of this Lie algebra which corresponds to
a  trivial representation summand). We
will use this same notation $Ad^{\Delta^A}$ when restricting to various Lie
subgroups of
$H_{\mathbb
R} L^A$ containing $L^A$. The notation $ Ad^{w(\Delta^A)}$  with $w\in W$ refers
to the similar construction with respect to $w(\Delta^A)$. When $A=\Pi$ and
$\Delta^A=\emptyset$ then $ Ad^{w(\Delta^A)}$ will refer to the (trivial)
one dimensional
representation of $\{e \}$.

\begin{defn}
\label{D4.1.1}
Let  $H_{\text{fund}}^A$ be defined by:
\begin{equation}
\exp\left( \left\{\   \underset { \alpha_i \not \in A }  {\sum } c_i m_{\alpha_i}^{\circ} \, : \ c_i\in {\mathbb R}  \right\} \right)=H^A_{\text{fund}}~.
\label{D4.1.3}
\end{equation}
Note that  usually ${H^A}\not=H_{\text{fund}}^A$.
\end{defn}
Then we have:

\begin{prop}
The images of $Ad^{\Delta^A}$ on subsets of Cartan subgroups satisfy

\begin{enumerate}
\item 
$Ad^{\Delta^A}({\mathcal D} (\Delta ))={\mathcal D} (\Delta ^A )$, \label{P4.1.1b}

\item  $Ad^{\Delta^A}(H^A_{\mathbb R})=H_{\mathbb R}(\Delta^A),$ \label{P4.1.1c}

\item $Ad^{\Delta^A}( H_{\text{fund}}^A )=Ad^{\Delta^A}(H^A)$. \label{P4.1.1d}
\end{enumerate}
\label{P4.1.1}
\end{prop}
\begin{Proof}
We first point out that the exponential map in the Lie groups $L^A_{\mathbb R}$,
$Ad^{\Delta^A}$  gives a diffeomorphism between ${\mathfrak h}^A$ and the
corresponding connected Lie group. Thus  the Lie groups $H^A$, $
Ad^{\Delta^A}(H^A)$  are isomorphic. This takes care of  part \ref{P4.1.1c})
on the
level of the connected component of the identity. 

Also the image under $ad^{\Delta^A}$ of the set $\{ m_{\alpha_i}^\circ: \alpha_i \not\in A \}$ gives rise to a basis of the Cartan subalgebra of $ad^{\Delta^A} ({\mathfrak  l}^A )$.  Thus exponentiating we obtain part \ref{P4.1.1d}).

Recall that  $\chi_{\alpha_i}(h_j)=\exp(\pi  \sqrt {-1}\delta_{i,j})$ and
$Ad^{\Delta^A}(h)=I$ if and
only if $\chi_{\alpha_i}(h)=1$ for all $\alpha_i\in \Pi\setminus A$. We have thus
$Ad^{\Delta^A}(h_i)=I$ for $\alpha_i\in A$ and $\{\ Ad^{\Delta^A}
(h_\epsilon ):
\epsilon \in {\mathcal E} \}= \{ Ad^{\Delta^A} (h_\epsilon ): \epsilon \in {\mathcal E}^A
\}={\mathcal D} (\Delta ^A )$.   The elements $h_i$ with $\alpha_i\in A$ are in the
center of $H_{\mathbb R} L^A$. This proves part \ref{P4.1.1b}).

For part \ref{P4.1.1c}), we proceed by  noting that by definition of $H_{\mathbb R} (\Delta^A) $ it must
be generated by  the Lie group $Ad^{\Delta ^A}(H^A)$ with Lie algebra
${\mathfrak h}^A$
and the elements $Ad(h_i)$ with $\alpha_i\in \Pi\setminus A$ (playing the role that the
$h_i$ play in the definition of $H_{\mathbb R}$). This is the same as 
$Ad^{\Delta ^A}(H^A_{\mathbb R})$.
\end{Proof}

\subsection{Action of the Weyl group on oriented colored Dynkin diagrams}

\begin{prop}
The action of $W_S$ on ${\mathbb D}(S)\times
\{\ \pm
1 \}$ in Definition \ref{D3.3.1} is well-defined. \label{P4.1.2}
\end{prop}

\begin{Proof}
Recall that $W$ is a Coxeter group, (Proposition 3.13  \cite{kac:90}) and it
thus has defining relations $s_{\alpha_i}^2=e$ and $(s_{\alpha_i}s_{\alpha_j})^{m_{ij}}=e$ where 
 $m_{ij}$ is $2,3,4,6$ depending on the number of lines joining $\alpha_i$ and $\alpha_j$
 in the Dynkin diagram.  The case $m_{ij}=2$ occurring when $\alpha_{i}$ and $\alpha_{j}$  are not connected in the Dynkin diagram. The only relevant cases then are when $\alpha_i$, $\alpha_j$ are both {\em colored} and connected in the Dynkin diagram. The only relevant vertices in the Dynkin diagram are those connected with these two and which are uncolored. We are thus reduced to very few  non-trivial possibilities: $D_5$, $A_4$, $F_4$, $B_4$, $C_4$; smaller rank cases being very easy cases.  We verify only one of these cases , the others being almost identical with no additional difficulties.  Consider the case of $D_5$ where  $\alpha_i=\alpha$, $\alpha_j=\beta$  are the two simple roots which are not \lq\lq endpoints\rq\rq   in the Dynkin diagram and all the others are uncolored. Then $r_{\alpha}=1$ but $r_{\beta}=2$. For book-keeping purposes it is convenient to temporarily represent {\em red} as $-1$ and {\em blue} as $1$ and simply follow what happens to these two roots and an orientation $o=1$. Hence all the information can be encoded in a triple $(\epsilon_1, \epsilon_2, o)$ representing the colors of these two roots and the orientation $o$. Now we apply $(s_\alpha s_\beta)^3$ to $(\epsilon_1,\epsilon_2,o)$:
\begin{equation}
\begin{array} {cc}
& \displaystyle{(\epsilon_1, \epsilon_2, o=1)\overset{s_{\beta}}{\longrightarrow}(\epsilon_1 \epsilon_2,\epsilon_2, (\epsilon_2)^2o=1)\overset{s_{\alpha}}
{\longrightarrow}
} \\
& \displaystyle{(\epsilon_1\epsilon_2,\epsilon_1, (\epsilon_1\epsilon_2)^1  o)\overset{s_{\beta}}{\longrightarrow}(\epsilon_2, \epsilon_1, (\epsilon_1)^2\epsilon_1\epsilon_2 o= \epsilon_1\epsilon_2 o)\overset{s_{\alpha}}
{\longrightarrow}}\\
& \displaystyle{(\epsilon_2, \epsilon_1\epsilon_2, (\epsilon_2)^1 \epsilon_1\epsilon_2 o=\epsilon_1o) \overset{s_{\beta}}{\longrightarrow} ( \epsilon_1, \epsilon_1\epsilon_2, (\epsilon_1\epsilon_2)^2\epsilon_1 o=\epsilon_1 o) \overset{s_{\alpha}}{\longrightarrow} }\\
& \displaystyle{(\epsilon_1, \epsilon_2, o) }
\end{array}
\nonumber
\end{equation}
\end{Proof}

\subsection{The set ${H}_{\mathbb R}^{\circ}$}
\label{ss4.2}  
By Proposition \ref{P4.1.1}, if $\chi^{\Delta^A}_{\alpha_i}$, $\alpha_i\in
\Pi\setminus A$ is
a  root character of  $H_{\mathbb R}(\Delta^A)$ acting on ${\mathfrak l}^A$ , then it is
related to $\chi_{\alpha_i}$ by $\chi^{\Delta^A}_{\alpha_i}\circ
Ad^{\Delta^A}=\chi_{\alpha_i}$. The map $\chi^{\Delta^A}_{\alpha_i}$, $\alpha_i\in\Pi\setminus A$ provides the local coordinates on $H_{\mathbb R}^{\circ}$ which consists of the Cartan subgroups $H_{\mathbb R}(\Delta^A)$.
\begin{defn}
 \label{D4.2.1} 
Let $H_{\mathbb R}^{\circ}$ be defined as 
\begin{equation}
 {H}_{\mathbb R}^{\circ}= \underset { A\subset
\Pi }    {\bigcup } H_{\mathbb R}  (\Delta^A) \times \{[e]^{A} \} \, .
\nonumber
\end{equation}
We then define a map $\phi_e : H_{\mathbb R}^{\circ}\to {\mathbb R}^l$
 as follows:
\begin{equation}
\phi_e(h, [e]^A)=(\phi_{e,1}(h),.., \phi_{e,l}(h)) ,
\nonumber
\end{equation}
where
$\phi_{e,i}(h)=\chi^{\Delta^A}_{\alpha_i}(h)$ whenever $\alpha_i \not \in
A$ and$\phi_{e,i}(h)=0$ if $\alpha_i \in A$. Denote $\phi_e^A$ the
restriction
of
$\phi_e$ to $H_{\mathbb R} (\Delta^A) \times \{ [e]^{A} \}.$
\end{defn}

By Proposition \ref{P4.1.1}  part c) we can compose with $Ad^{\Delta^{A}}\times 1$
and  re-write the  domain of $\phi_e\circ (Ad^{\Delta^A}\times 1)$   as $
\underset { A\subset \Pi } {\bigcup } H_{\mathbb R}^A \times \{ [e]^{A}
\} \,$.  We also define
\begin{equation}
 \phi_w= (\phi_{w,1},..,\phi_{w,l})~:~w(H_{\mathbb R}^{\circ})\longrightarrow
 {\mathbb R}^l~,
 \nonumber
\end{equation}
with
\begin{equation}
 w({H}_{\mathbb R}^{\circ} ):=\underset { A\subset\Pi }
{\bigcup } H_{\mathbb R} (w(\Delta^A)) \times \{ [w]^{A} \} \ , 
\nonumber
\end{equation}
by setting
\begin{equation}
\phi_{w,i}(x, [w]^A)=\chi^{w(\Delta^A)}_{w\alpha_i}(Ad^{w(\Delta^A)}(wh))=
\chi_{w\alpha_i} (wh)=\chi_{\alpha_i}(h),
\nonumber
\end{equation}
 where $x=Ad^{w(\Delta^A)}(wh)\in H_{\mathbb R}
(w(\Delta^A) )$ with $h\in H_{\mathbb R}$, if
$\alpha_i\not\in A$. For the case when $\alpha_i\in
A$, we set 
\begin{equation}
\phi_{w,i}(x,[w]^A)=0.
\nonumber
\end{equation}

\subsection{The sets $H^A_{\mathbb R}$}
 \label{ss4.3}
 We here note an isomorphism between several presentations  of  a split Cartan subgroup
 of a Levi factor.

\begin{prop}
All the following are isomorphic as Lie groups:

\begin{enumerate}


\item  $H^A_{\mathbb R}$, \label{P4.3.1c}

\item  $ Ad^{\Delta^A} (H_{\mathbb R} )$, \label{P4.3.1d}

\item   $Ad^{\Delta^A} (H^A_{\mathbb R})$, \label{P4.3.1e}

\item  $H_{\mathbb R} (\Delta^A)$. \label{P4.3.1f}
\end{enumerate}
\label{P4.3.1}
\end{prop}

\begin{Proof} 
The groups
in \ref{P4.3.1e} and \ref{P4.3.1f} are isomorphic by  proposition \ref{P4.1.1}. We can write ${\mathfrak h}_{\mathbb C}=
{\mathfrak z} + {\mathfrak h}^A_{\mathbb C}$ where  ${\mathfrak z}$ is defined by  $z\in {\mathfrak z}$
if and
only if $\langle \alpha_i , z\rangle =0 $ for all  $\alpha_i\in \Pi\setminus A$.

Now the group  $\exp({\mathfrak z})$ is the center   of $H_{\mathbb C} L^A_{\mathbb C}$.
Intersecting with $\tilde G$  we obtain the center of $H_{\mathbb R} L^A$. To find
this intersection with $\tilde G$, we must find all  $x+\sqrt {-1} y\in {\mathfrak z}$
such that
$e^{\langle \alpha_j, x+\sqrt {-1} y\rangle }$ is real for all  $\alpha_i \in \Pi$ (see the
proof of
Lemma \ref{L2.2.1}).  We find that $x\in {\mathfrak z} \cap
{\mathfrak h}$ and that $y$ is an integral linear combination of the $y_i$. As in the
proof of Proposition \ref{P4.1.1} the elements $\exp(\sqrt {-1} y_i)$ which are in the center
of
$H_{\mathbb R } L^A$ are those for which $\alpha_i \in A$.  The center of  $H_{\mathbb R}
L^A$ is then the group generated by $h_i$ with $\alpha_i \in A$ and $\exp({\mathfrak
z} \cap {\mathfrak h})$. Since $\exp({\mathfrak z} \cap {\mathfrak h}) \cap H^A=\{ e \}$, we
have that
 $H_{\mathbb R}$ divided by $\exp({\mathfrak z} \cap {\mathfrak h})$ is isomorphic to ${\mathcal D}{
H}^A_{\mathbb R}$. The  groups $H^A_{\mathbb R}$,   ${\mathcal D}{
H}^A_{\mathbb R}$ and those involved in  \ref{P4.3.1d}, \ref{P4.3.1e} or \ref{P4.3.1f}  all differ by a subgroup of the ${\mathcal D}\subset  \exp({\mathfrak z}) $ which is annihilated by $Ad^{\Delta^A}$. From here and Proposition \ref{P4.1.1} part \ref{P4.1.1b}),
the isomorphism betwen   \ref{P4.3.1d}, \ref{P4.3.1e}  and  \ref{P4.3.1f} follows.
\end{Proof}

\section{The set ${\hat H}_{\mathbb R}$}
 \label{S5} 
 We here define our main object ${\hat H}_{\mathbb R}$ as a union of Cartan subgroups of Levi factors and their $W$-translations.

Let $\hat W$ be the disjoint union of all the quotients $W/W_{\Pi\setminus A}$ over $A\subset \Pi$. Each of the elements $[w]^A\in W/W_{\Pi\setminus A}$ parametrizes a
parabolic subgroup. First $[e]^A$ corresponds to a standard parabolic
subgroup of
$G_{\mathbb C}$.  (Proposition  7.76 or Proposition 5.90 of \cite{KN}). This is just the
parabolic subgroup determined by the subset $\Pi\setminus A$ of $\Pi$. Then we translate
such a parabolic subgroup with $w\in W$. The resulting parabolic subgroup
corresponds to $[w]^A$. Thus $\hat W$ is the set of all the $W$-translations of
standard parabolic subgroups.
\begin{defn}
\label{D5.1.1}
We define
\begin{equation}
\hat H_{\mathbb R}=\displaystyle{\bigcup_{A\subset \Pi}
\bigcup_{w\in W} H_{\mathbb R} (w(\Delta^A)) \times \{[w]^A \} \,}.
\label{5.1.1}
\end{equation}
\end{defn}
From Proposition \ref{P4.3.1}, $ H^A_{\mathbb R}$ can be replaced by $H_{\mathbb R} (\Delta^A)$ or $Ad^{w(\Delta^A)}(H_{\mathbb R})$.  Thus we can alternatively  write $\hat H_{\mathbb R}$ as a subset of $H_{\mathbb R}\times \hat
W$. We recall $w^{\bullet,A}$ or just $w^\bullet$  as in Notation
\ref{N3.3.1}.
We then have:
\begin{equation}
\hat H_{\mathbb R}=\displaystyle{\bigcup_{A\subset \Pi}
\bigcup_{w^{\bullet,A}\in W/W_{\Pi\setminus A}}  w^{\bullet,A}
(H^A_{\mathbb R}) \times \{[w]^A \} \,}. 
\label{5.1.2}
\end{equation}

Also fixing an isomorphism $\xi_w:H_{\mathbb R} (\Delta^A)\to H_{\mathbb R} (w(\Delta^A))$
inducing a set bijection (by composition) $\xi_w^*:w(\Pi)\to \Pi$, we have
\begin{equation}
{\hat H}_{\mathbb R}\cong\underset {A\subset \Pi} { \bigcup } W \underset {W_{\Pi\setminus A}} {\times }
H_{\mathbb R}(\Delta^A) ~,
\nonumber
\end{equation}
where
an element $(w^{\bullet, A}, h)$ on the right hand side is sent to $(\xi_w(h), [w]^A)$. This endows $\hat H_{\mathbb R}$ with a $W$ action.

In what follows it is useful to think of a colored Dynkin diagram (e.g.,
$\circ_R-\circ$) as parametrizing a \lq\lq box\rq\rq (e.g., $[-1,1]$) We
need to
further subdivide this box into $2^l$ smaller boxes by dividing it into regions
according to the sign of each of the coordinates (e.g., $[-1,0]$ $[0,1]$). We
also need to consider the boundary between these $2^l$ regions (e.g., $\{ 0
\}$). We will introduce an additional sign or a zero to keep track of such
subdivisions (e.g., $\circ_R - \circ_+$, $\circ_R-\circ_-$, $\circ_R-\circ_0$
for $[-1,0], [0,1], \{0\}$) (see also Example in Section \ref{SSE5.2.2}).
We will then do the same thing with a colored Dynkin diagram of the form
$(D,[w]^{\Pi\setminus S})$ by assigning labels in $\{\pm 1,0 \}$ to the vertices in
$\Pi\setminus S$.

The main purpose of the following is to associate certain sets in ${\hat H}_{\mathbb R}$ 
to the colored Dynkin diagrams, the signed-colored Dynkin diagrams and the sets associated 
to them only play an auxiliary role.

We introduce the following notation.  This notation is illustrated in Example
in Section \ref{SSE5.2.2} and Figure \ref{F5.2.1}.

\begin{notation}
\label{D5.1.2}

A {\em signed-colored Dynkin diagram} $\check D$ is a Dynkin diagram
with some
vertices colored ($R$ or $B$) and the remaining vertices labled $+,-$ or $0$.
The followings are auxiliary objects to keep track of signs, zeros and colors. 
\begin{itemize}

\item  $\check \eta :\Pi\to \{ \pm 1, 0 \}$ function which agrees
with $\eta$ on $S$ and determines the sign labels in  $\check D$

\item  $A=A(\check D)=A(\check\eta)=\{ \alpha_i\in \Pi: \check\eta
(\alpha_i)=0 \}$ 

\item $K(\eta)$ the set of all $\check \eta:\Pi\to \{
\pm 1, 0\}$ which agree with $\eta$ in $S$

\item $\epsilon_{\check\eta}
\in {\mathcal E}^A$ the element which agrees with $\check\eta$ on $\Pi\setminus A$ ($A=A(\check D)$)

\item  $\check D
(S,A,\epsilon_{\check \eta}, [w]^{\Pi\setminus S})$ the unique signed-colored Dynkin diagram
 attached to a colored Dynkin diagram $(D,[w]^{\Pi\setminus S})$ or to
$(S,\eta,[w]^{\Pi\setminus S})$,
where the vertices in $\Pi\setminus S$ are given a label in $\{ \pm 1,0 \}$.  
\end{itemize}

We now associate a subset of $\hat H_{\mathbb R}$ to a signed-colored Dynkin diagram
with $S=\emptyset$. Notice that in this case $\epsilon_{\check \eta}$ can
be any
element of ${\mathcal E}^A$. Recall that $\epsilon_{\check\eta}(\Delta^A)$ is
then the
element of ${\mathcal E}(\Delta^A)$ that corresponds.

We associate to a signed colored Dynkin diagram $(\emptyset, A, \epsilon_{\check\eta} ,[e]^{\Pi\setminus S})$  two sets, one includes walls and the other doesn't (in order to avoid duplicate notation we denote the signed-colored Dynkin diagram and the set with the same notation): 
$$ {   (\emptyset, A, \epsilon_{\check\eta} ,[e]^{\Pi\setminus S})}^{\le}
= H(\Delta^A)_{
\epsilon_{\check\eta}(\Delta^A) }^{\leq}  \times \{[e]^A \}  $$
When $A=\emptyset$ (no zeros), these are the $2^l$ boxes in the antidominant chamber.
We define the second related set as:

$$ {   (\emptyset, A, \epsilon_{\check\eta} ,[e]^{\Pi\setminus S})}
= H(\Delta^A)_{
\epsilon_{\check\eta}(\Delta^A) }^{<}  \times \{ [e]^A \}  $$

The chamber walls of the antidominant chamber of the Cartan subgroup are defined as:

\begin{itemize}

\item  ${ D} (\alpha_i,A,\epsilon)^{\leq }=    \{\ h\in
H(\Delta^A)^{\leq}_{\epsilon (\Delta^A)} : |\chi_{\alpha_i}(h)|=1 \}$
(the $\alpha_i$-wall)

\item  $ D(\alpha_i,A,\epsilon )^{<}= { D} (\alpha_i,A,\epsilon)^{\leq}\cap  \{\ h\in
H(\Delta^A)^{\leq}_{\epsilon (\Delta^A)} : |\chi_{\alpha_j}(h)| < 1 \text { if } j\not=i,\alpha_j \in \Pi\setminus A \}$

\end{itemize}

 We next consider the case
of  $S=\{\alpha_i \}\not\in A$ and then the general case of any $S\subset
\Pi$
with $A\subset \Pi\setminus S$ and any $\epsilon_{\check\eta}\in {\mathcal E}^A$. This defines the walls for Levi factor pieces corresponding to subsystems of the Toda lattice (we here list open walls):

\begin{itemize}

\item  $  (\{\ \alpha_i \},A,\epsilon_{\check\eta}) = D
(\alpha_i,A,\epsilon_{\check\eta} )^{<}  \times \{ [e]^A \}$.

\item $(S,A,\epsilon_{\check \eta})=\underset
{\alpha_i\in S}{\bigcap} (\{\ \alpha_i
\},A,\epsilon_{\check\eta})$.

\end{itemize}

We now associate a set in ${\hat H}_{\mathbb R}$ to a colored Dynkin diagram. We define a set denoted
$D=(S,\epsilon_\eta)$ as follows:
For $S\not=\emptyset$ (so that there is an $\eta:S\to \{ \pm 1\}$),

\begin{itemize}

\item $(S,\epsilon_\eta)=\underset
{\check\eta\in K(\eta)} {\bigcup}(S,A(\check\eta),\epsilon_{\check \eta})$, and
if $S=\emptyset$, 
$(\emptyset ,\epsilon_o)=\underset
{{\check\eta}\in \{\pm 1,0\}^l}{\bigcup} (\emptyset, {\epsilon_{\check\eta}}) $
\end{itemize}
Here $\epsilon_o=(1,\cdots,1)$, and see (\ref{5.2.4}) and (\ref{5.2.20}) in Section\ref{SSE5.2.2}.

We consider the $w$-translations of the colored Dynkin diagrams: For this
write $w$ uniquely as $w=w^\bullet w_\bullet$ with $w_\bullet \in W_S$.
\begin{itemize}

\item  $(w_\bullet D,[w]^{\Pi\setminus S})=w((D,[e]^{\Pi\setminus S})). $

\end{itemize}

In order to define the $W$-translations of sets associated to signed-colored Dynkin diagrams we have to extend the definition of the $W_S$-action to the
signed-colored Dynkin diagrams. The definition is exactly the same if we treat
the label $-$ as if it were an $R$ and the label $+$ as if it were a $B$ as in
Definition \ref{D3.3.1}. Then consider $(w_\bullet \check D,[w]^{\Pi\setminus S})$ and let:

\begin{itemize}

\item  $(w_\bullet \check D,[w]^{\Pi\setminus S})=
w( (\check D,[e]^{\Pi\setminus S})). $
\smallskip

\end{itemize}

\end{notation}

\begin{rem}\label{ss6.1}

 We refer to the set
 $\overline {(\emptyset,\epsilon_o) } \subset \hat H_{\mathbb R}$ as
the $antidominant$ $chamber$ of ${\hat H}_{\mathbb R}$. 
Recall Notation \ref{N3.3.1} $w^{\bullet,A} \in [w]^A$. The following justifies
our definition of the \lq\lq antidominant chamber \rq\rq of ${\hat H}_{\mathbb R}$ using  (\ref{5.1.2}) in the definition of ${\hat H}_{\mathbb R}$.
\end{rem}
\begin{prop}
We have
\begin{equation}
{\hat H}_{\mathbb R} = \underset { A\subset \Pi}{\bigcup}
\underset{w^{\bullet , A} \in W \atop
 \sigma \in W_{\Pi\setminus A} }     { \bigcup } w^{\bullet, A}
\left( \sigma \left(  H_{\epsilon (\sigma)}^{A,\leq }\right)\right) \times \{ [w]^A \} \, . 
\end{equation}
\label{P6.1.1}
\end{prop}

\begin{Proof} We set $w^\bullet= w^{\bullet, A}$ for simplicity. By  Proposition
\ref{P2.2.4} part b), we have that each $H_{\mathbb R}^A$  can be written as a union over
$\sigma\in W_{\Pi\setminus A}$ of sets of the form $\sigma \left(H_{\epsilon (\sigma) }^{A,\leq }\right)$
and thus by  the definition (\ref{5.1.2}) of ${\hat H}_{\mathbb R}$, we
conclude the statement.
\end{Proof}

\section{Colored Dynkin diagrams and the corresponding cells}
 \label{S7}
Here we consider the manifold structure and the topology of ${\hat H}_{\mathbb R}$ as the union of cells parametrized by colored Dynkin diagrams.
\subsection{Action of the Weyl group on the sets $(S,\epsilon_{\eta})$}
 \label{ss7.1}
 Let ${\mathcal M}^{geo}$ be the complex with the sets $(D, [w]^{\Pi\setminus S})$. We then consider the action of $W_S$ on the union of all the sets $D=(S,\epsilon_{\eta})$ having a fixed nonempty set $S$ of colored vertices and endowed with an orientation $o$.
Similarly  we can endow each of the terms in the chain complex ${\mathcal M}^{geo}$ with a action of $W$ by translating the sets corresponding to colored Dynkin diagrams with the $W$-action and taking into account changes of orientation induced on the oriented boxes. This new action  of $W$ on ${\mathcal M}^{geo}$ could in principle be different from the $W$-action on the chain complex ${\mathcal M}_*$ (see Proposition \ref{P3.3.2}).

We now become more explicit about the $W$ action that was just introduced on
$\mathcal M^{geo}_*$:
Note that since
$\chi_{\alpha_{j_i}}(h)=\pm 1$  for any $\alpha_{j_i}\in S$ and $h\in \overline { (S, \epsilon_\eta)}$ (Notation \ref{D5.1.2}),
if $h\in   (S,\epsilon_\eta)$ then $s_{\alpha_{j_i} } \chi_{\alpha_j}
(h)=\chi_{\alpha_j}(h)\chi_{\alpha _{j_i}} (h)^{-C_{j,j_i}}= \pm
\chi_{\alpha_j}(h)$.  Hence, in terms of the coordinates given by
$\phi_e$, the
action of $s_{\alpha_{j_i}}$ is given by a diagonal matrix whose non-zero
entries
are $\pm 1$. This matrix has some entries corresponding to the set $S$ and
other
entries corresponding to $\Pi\setminus S$.  The entries corresponding to the set $S$
change by a sign as described in  Proposition  \ref{P2.2.3}. The statement in
Proposition  \ref{P2.2.3} just means that the set $ (S, \epsilon_\eta)$ is sent to
$ (S, \epsilon_{\eta^\prime })$ with $\eta^\prime=(\eta_{j_1},..,\eta_{j_s})$,  and
$\eta_{j_i}^\prime = \eta_{j_i}(-1)^{C_{j,j_i}}$.

The determinant of the diagonal submatrix corresponding to elements in
$\Pi\setminus S$ is
the sign $(\eta_i)^r$ with (see Definition \ref{D3.3.1})
\begin{equation}
r=\Big|\{\ \alpha_j\in \Pi\setminus S \ : \ C_{j,i}
\text { is odd } \}\Big|.
\nonumber
\end{equation}

Consider the {\em set} $ (S,\epsilon_\eta) $ endowed with a fixed orientation $\omega$
corresponding to $o=1$. Then $s_{\alpha_i}$ with $\alpha_i\in S$ sends $(S,\epsilon_\eta) $ to $ (S,\epsilon_{\eta^\prime}) $ endowed with the orientation  $(\eta_i)^r
\omega$. We now consider the $\mathbb Z$ module of  formal integral
combinations of
the sets $ (S,\epsilon_\eta)$,
\begin{equation}
{\mathbb Z}\left[\ (S,\epsilon_\eta)\ :~ S\subset \Pi,~|S|=k,~\epsilon_{\eta}\in{\mathcal E}^{\Pi\setminus S}~\right].
\nonumber
\end{equation}

We keep track of the orientation by putting a sign $\pm$ in front of $
(S,\epsilon_\eta)$. This $\mathbb Z$ module acquires a ${\mathbb Z}[W_S]$ action that
corresponds
to the abstract construction given  in Definition \ref{D3.3.2} with colored Dynkin
diagrams. By considering all the $W$-translations of the $ (S,\epsilon_\eta)$ we
generate the module denoted above by ${\mathcal M} (S)$. The direct sum of all these
$\mathcal M(S)$ over $|S|=k$ is denoted $ {\mathcal M}_{l-k}$ (see Definition \ref{D3.3.2}). Hence there is no difference as $W$-modules between
${\mathcal M}^{geo}_*$ and ${\mathcal M}_*$.

We can now summarize this discussion in the following:
\begin{prop} The action of $W$ on ${\mathcal M}_*$ (subsection \ref{ss3.3}) and the action 
of $W$ on  ${\mathcal M}^{geo}$ are isomorphic. 
\label{P7.1.1}
\end{prop}
We will then drop the superscript ${\cdots }^{geo}$ from the notation in view of Proposition \ref{P7.1.1}.

\subsection {Manifold structure on ${\hat H}_{\mathbb R}$}
 \label{ss7.2}

Recall the  map $\phi_e$ in Definition \ref{D4.2.1} whose domain is $\overset
{\circ} {H}_{\mathbb
R}= \underset { A\subset \Pi }  {\bigcup } H_{\mathbb R} (\Delta^A) \times \{
[e]^A \} \,$  and co-domain is ${\mathbb R}^l$.
We also  have defined $\phi_w$ , with domain $ \underset { A\subset\Pi}
 {\bigcup } H_{\mathbb R} ( w(\Delta^A)) \times \{[w]^A \} \,) $.
 We will  use these
maps to give ${\hat H}_{\mathbb R}$ coordinate charts leading to a manifold structure. We then have the following three Propositions:

\begin{prop}
The image  $\phi_e (H_{\mathbb R}
(\Delta^A))\times \{
[e]^A \}$ consists of all $(t_1,..,t_l)\in { \mathbb R}^l $ such that $t_i  \not=
0$ if
and only if $\alpha_i \not \in A$. The map $\phi_e$ is a bijection between $
\overset {\circ}  {H}_{\mathbb R}=\underset { A\subset\Pi }   {\bigcup }
(H_{\mathbb R} (\Delta^A) ) \times \{ [e]^A \} \,$ and ${\mathbb R}^l$.  \label{P7.2.1}
\end{prop}

\begin{Proof}
 We start with  the last statement, that $\phi_e$ is a
bijection; we
have that $\phi_e$ is injective because the scalars
$\chi^{\Delta^A}_{\alpha_{j_1}}(h),..,\chi^{\Delta^A}_{\alpha_{j_m}}(h)$
determine all the root
characters $\chi^{\Delta^A}_\phi (h)$, $\phi\in \Delta^A$ for $h\in H_{\mathbb R}
(\Delta^A)$ and
these scalars determine $h$ in the adjoint group (Remark \ref{R2.2.1}). From
Proposition \ref{P4.1.1}  part d)  it follows that we can regard $H_{\mathbb R} (\Delta^A)$ as
$Ad^{\Delta^A}(H_{\mathbb R})$ (see (\ref{D4.1.3})). We prove surjectivity by
proving first the
statement
concerning the image $\phi_e\circ (Ad^{\Delta^A}\times 1) ( H_A \times \{
[e]^A
\})$ . When all the sets $A$ are considered then all of
$\mathbb
R^l$  will be seen to be in the image of $\phi_e$.   First consider
$h_\epsilon
h$ with $h= \exp(  \underset { \alpha_i \not \in A } { \sum } c_i
m_{\alpha_i}^\circ \,)  $ and $\epsilon\in \mathcal E^A$.   We now apply $\phi_e
\circ (Ad^{\Delta^A}\times 1 ) $. Since $\langle \alpha_i , \underset {
\alpha_j \not
\in A }   { \sum } c_j m_{\alpha_j}^\circ \,\rangle=  c_i \frac
{(\alpha_i,\alpha_i)}{2} \, $ we obtain, by exponentiating,  $\chi_{\alpha_i}
(h_\epsilon h)=\epsilon_i e^{c_i\frac {(\alpha_i,\alpha_i)}{2}}$.   The set
$\phi_e (H_{\mathbb R} (\Delta^A)\times \{[e]^A \}) $ becomes the image of the map:
${\mathbb R}^l \to {\mathbb R}^l$ given by first defining a map  that sends
$(\epsilon_1t_1,..,\epsilon_l t_l)\to (f_1,..,f_l)$ with $f_i=\epsilon_i
t_i^{\frac {(\alpha_i,\alpha_i)}{2}}$  for  $t_i>0$. This map is modified  so
that whenever $\alpha_i\in A$ then the $i$-th coordinate is replaced with
$0$. We
denote this modified map by  $F^A$.  The domain and the image of $F^A$
therefore
consists of the set $\{\ (s_1,..,s_l): s_i=0 \text { if } \alpha_i\in A \}$.

Together all these $F^A$ give rise to one single map $F:{\mathbb R}^l \to {\mathbb R}^l$
which is surjective.
\end{Proof}

\begin{prop}
The image  $\phi_w ( H_{\mathbb R}
(w(\Delta^A))\times
\{ [w]^A \})$ consists of all $ (t_1,..,t_l)\in {\mathbb R}^l $ such that $t_i
\not=
0$ if and only if $\alpha_i \not \in A$. \label{P7.2.2}
\end{prop}

\begin{Proof} This follows from Proposition \ref{P6.1.1} and the fact that
$\phi_{e,i}(w(h),A)=\chi_{w\alpha_i} (w(h))= \chi_{\alpha_i} (h)$ for
$\alpha_i\not\in A$.
\end{Proof}
\begin{prop}
The image  $\phi_w ((\emptyset,\epsilon_o))$
 consists of all $ (t_1,..,t_l)\in{ \mathbb R}^l $ such that $-1 < t_i < 1$.
 The sets $\phi_w ( (S,\epsilon_\eta, [w]^{\Pi\setminus S}))$ as $S\subset
\Pi$, $S\not=\emptyset$
 varies, give a cell decomposition of the boundary of the box $[-1,1]^l$.
 In particular, the sets $(S,\epsilon_\eta, [w]^{\Pi\setminus S})$ give a cell
decomposition
 of the smooth manifold $\hat H_{\mathbb R}$. \label{P7.2.3}
\end{prop}

\begin {Proof} This follows from Proposition \ref{P7.2.2} but is better understood in 
Example  \ref{E5.2.2}. We omit details.
\end{Proof}

\begin{rem}
\label{R7.2.1}
There is a more convenient cell decomposition of ${\hat H}_{\mathbb R}$
for the purpose of calculating homology explicitly. The only change is that
the $l$ dimensional cell becomes the union of all the $l$-cells together with
all the (internal) boundaries corresponding to colored Dynkin diagrams
where all the colored
vertices are colored $B$. This is the set:
\begin{equation}
{\hat H}_{\mathbb R} \setminus \underset
 {S\subset \Pi,~ w\in W \atop \eta \text { such that } \eta(\alpha_i)=-1
 \text { for some }\alpha_i\in S}
 {\bigcup } (S,\epsilon_\eta,[w]^{\Pi\setminus S}).
 \nonumber
\end{equation}
 This set can be seen
to be homeomorphic to
 ${\mathbb R}^l$. With this cell decomposition there is exactly one  $l$ cell;
 and the other lower dimensional cells correspond to colored Dynkin diagrams
 which are parametrized by pairs $(D,[w]^{\Pi\setminus S})$, such that, at least one
vertex
of $D$
 has been colored $R$. In terms of the $(S,\epsilon_\eta,[w]^{\Pi\setminus S})$,
 the top cell would instead be defined to consist of  ${H}_{\mathbb R}^{\circ}$.
(See Figure \ref{F3.3.1} for this remark.)

The big cell in this decomposition with a fixed orientation  corresponds to
the element
 $c_l=\underset { w\in W}{\sum} (-1)^{l(w)} (D,\epsilon_o, w)$
($S=\emptyset$).
 This element satisfies $\partial_l(c_l)=2(c_{l-1})$ for some $c_{l-1}\not=0$
 (except in the case of type $A_1$). Note that $(D,w)$ and $(D,ws_{\alpha_i}),$
 with $S=\emptyset$ and $l(ws_{\alpha_i})=l(w)+1$, appear with opposite signs
 in $c_l$. When $\partial_l$ is applied and  the $\alpha_i$ is colored $R$
 the sign $(-1)^{r_{\alpha_i}}$ in Definition \ref{D3.3.1} makes these two
 terms contribute as $2(D^\prime,[w]^{\Pi\setminus \{\ \alpha_i \}})$ with $D^\prime$
 the new colored Dynkin diagram obtained. The terms from the boundary
$\partial_l$
 obtained by coloring $\alpha_i$ with  $B$ will  cancel since the action of
$s_{\alpha_i}$
 on colored Dynkin diagrams
 is trivial when $\alpha_i$ is colored $B$. The case of $A_1$ is
 an exception because, in that case,  once $\alpha_1$ is colored $R$ no more
 uncolored vertices remain. The set $\{\ \alpha_j \in \Pi\setminus {\alpha_i}
:C_{j,i} \text { is odd } \}$
 is empty and $r_{\alpha_i}=0$. Thus there is cancellation in this case.
\end{rem}

\subsection {Topology on $\hat H_{\mathbb R}$ , coordinate charts, integral
homology}
  \label{ss7.3}

We define a topology on $\hat H_{\mathbb R}$ in which $U \subset  \underset { A\subset\Pi }  { \bigcup } H(w(\Delta^A)) \times \{ [w]^A \} \, $  is
open if
and only if $\phi_w (U)$ is open in the usual topology of ${\mathbb R}^l$. The maps
$\phi_w$ become coordinate charts and since the compositions $\phi_w \circ
\phi_\sigma ^{-1}$ are $C^\infty$ on their domain, then ${\hat H}_{\mathbb R}$
acquires
the structure of a smooth manifold.  The $W$ action becomes a smooth action.

\begin{defn}
\label{D7.3.1}
{\it Filtration of ${\hat H}_{\mathbb R}$ and the chain complex ${\mathcal M}^{CW}_*$}:
We construct a filtration of the topological
space ${\hat H}_{\mathbb R}$ in the sense of \cite{MU} p. 222.  Let $ X_{l-k}$ denote the
union of all the sets of the form $(D,[w]^{\Pi\setminus S})$ over all $w\in W$ and
$S\in P(\Pi)$ such that $|S|\ge k$. This is a closed set and $X_r\setminus X_{r-1}$ is a
union of sets of the form $(D,[w]^{\Pi\setminus S})$ with $|S|=l-r$.
The filtration $X_r$
$r=0,1,..,l$ satisfies the conditions of Theorem 39.4 in \cite{MU}.  We define a
chain
complex ${\mathcal M}^{CW}_*$ with boundary operators as in \cite{MU} 
\begin{equation}
\partial_r:
H_r(X_r, X_{r-1},{\mathbb Z}) \to H_{r-1}(X_{r-1}, X_{r-2},{\mathbb Z} ).
\nonumber
\end{equation}
\end{defn}

\begin{prop}
The smooth manifold ${\hat H}_{\mathbb R}$ is
compact, non-orientable (except if $\mathfrak g$ is of tupe $A_1$).
 The homology of the chain complex ${\mathcal
M}_*$, $H^k({\mathcal M}_*)= Ker \partial_k/image (\partial_{k-1})$ is isomorphic as a
${\mathbb Z}[W]$ module to  $H_k({\hat H}_{\mathbb R},{\mathbb Z})$. \label{P7.3.1}
\end{prop}

 \begin {Proof} The manifold ${\hat H}_{\mathbb R}$ is the finite union
of the
chambers as in Proposition \ref{P6.1.1}. Since the $W$ action on ${\hat H}_{\mathbb R}$
is by
continuous transformations, it then suffices to observe that the antidominant
chamber  is compact.  The antidominant chamber  $\overline {(\emptyset,\epsilon_o)}$ of 
Definition \ref{D5.1.2} can be seen to be compact by
describing explicitly its  image under $\phi_e$.  This image is a
\lq\lq box\rq\rq inside ${\mathbb R}^l$, as can be seen in Propositions \ref{P7.2.1},
\ref{P7.2.3}
namely
the set $\{\ (t_1,..,t_l): -1\leq t_i \leq 1 \}$.
The space ${\hat H}_{\mathbb R}$
is now
the finite union of the $W$ translates of this compact set. That the
boundary operators of ${\mathcal M}^{CW}$ agree with the boundary operators of
${\mathcal M}_*$
will follow from the fact that in the $\phi_w$
coordinates
the $(D,[w])$  is a \lq\lq box\rq\rq which is itself  part of the
boundary
of a bigger \lq\lq box\rq\rq (Proposition \ref{P7.2.3}). We start with the set
$\{\ (t_1,..,t_l): -1\leq t_i \leq 1 \}$ and note that
its boundary is combinatorially described by (\ref{3.2.4}) or
(\ref{3.3.6}). Note that $(D,w)$ (with
$S=\emptyset, [w]^{\Pi}=w$) represents the open box.
The faces are parametrized by coloring each of the
$l$ vertices $R$ or $B$ which then represent {\it opposite}
faces in the boundary. The signs are just chosen so that
$\partial_{k-1}\circ\partial_k=0$
for $k=1,\cdots,l$. This description may be best understood
by working out Example \ref{E5.2.2}.

All the cells thus
appear by taking the faces of a box $[-1,1]^l$ and then faces of faces etc.
By the same process of coloring uncolored vertices
$R$ or $B$ which give rise, each time, to a pair of
opposite faces in a box.
In each case (\ref{3.2.4}) or (\ref{3.3.6}) correctly describe the process of taking
the boundary of a box. Note that it is enough to study
what happens when $w=e$ and then consider the $W$ translates.

We now use Theorem 39.4 of \cite{MU} to conclude that
${\mathcal M}^{CW}_*$ computes integral homology. However each $\mathbb Z[W]$-module
appearing in  ${\mathcal M}^{CW}_*$ in a fixed degree, can easily be seen to
be identical with the corresponding term in ${\mathcal M}_*$.
By  Proposition \ref{P7.1.1} and the agreement of the boundary
operators, we obtain that
${\mathcal M}_*$
computes integral homology.  The non-orientability follows
if we use the second cell decomposition described in Remark \ref{R7.2.1} b) and c).
The unique top cell then has a non-zero boundary (except in the case of
${\mathfrak g}={\mathfrak {sl}}(2, {\mathbb R})$).
\end{Proof}

\section{Toda lattice and the manifold ${\hat H}_{\mathbb R}$}
 \label{S8}
We now associate the nmanifold ${\hat H}_{\mathbb R}$ with the Toda lattice by
 extending the results of Kostant in \cite{KOS}.  We start with the definition of
  the variety $Z_{\mathbb R}$ of Jacobi elements on ${\mathfrak g}$ where the Toda lattice is defined.
\subsection {The variety $Z_{\mathbb R}$ and
isotropy group ${\tilde G}^z$}
 \label{ss8.1}

\begin{defn}
\label{D8.1.1}
{\it Varieties of Jacobi elements}:
Let $S({\mathfrak g})$  be the symmetric
algebra of $\mathfrak g$. We may regard $S({\mathfrak g})$ as the algebra of polynomial
functions on the dual ${\mathfrak g}^\prime$.

If we consider the algebra of $G$-invariants of $S({\mathfrak g})$ , then by
Chevalley's theorem there are homogeneous polynomials $I_1,..,I_l$
in $S({\mathfrak g})^G$ which are algebraically independent and which generate
 $S({\mathfrak g})^G$. Thus
$S({\mathfrak g})^G$ can be expressed as ${\mathbb R}[I_1,..,I_n]$.

For $F={\mathbb C}$  or $F={\mathbb R}$, we consider the variety $Z_F$ of normalized
Jacobi
elements of ${\mathfrak g}_F$. Our notation, however, is slightly different from the
notation of \cite{KOS} in the roles of $e_{\alpha_i}$ and $e_{-\alpha_i}$.
Thus we let
\begin{equation}
\begin{array}{ll}
&J_F=\left\{\ X=x+ \displaystyle{ \sum_{i=1}^l
} (b_i
e_{-\alpha _i } +  e_{\alpha _i})\, : \, x\in {\mathfrak h},b_i\in F\setminus
\{0\}\right\},\\
&Z_F=\left\{\ X=x+ \displaystyle { \sum_{i=1}^l } (b_i
e_{-\alpha _i } +  e_{\alpha _i})\, : \, x\in {\mathfrak h},b_i\in F\setminus
\{0\}, X \in S(F) \right\} .\\
\end{array}
\nonumber
\end{equation}
We also allow {\em subsystems} which correspond to the cases
having some $b_i=0$:
\begin{equation}
\begin{array}{ll}
& \overset {\circ}  {J_F}=\left\{\ X=x+
\displaystyle { \sum_{i=1}^l } (b_i e_{-\alpha _i } +  e_{\alpha
_i})\, : \,x\in {\mathfrak h},b_i\in F \right\} , \\
& \overset {\circ} {Z_F}=\left\{\ X=x+
\displaystyle { \sum_{i=1}^l } (b_i e_{-\alpha _i } +  e_{\alpha
_i})\, : \, x\in {\mathfrak h},b_i\in F, X \in S(F) \right\} .
\end{array}
\nonumber
\end{equation}

\end{defn}

Kostant defines  in \cite{KOS} p. 218  a real  manifold $Z$ by considering all
elements $ x+ \displaystyle{ \sum_{i=1}^l } (b_i e_{-\alpha _i } +
e_{\alpha _i})$ in $Z_{\mathbb R}$  which in addition satisfy $b_i > 0$. We are
departing in a crucial way from  \cite{KOS} by allowing the $b_i$ to be
negative or even zero when $A\not=\emptyset$.
This extension gives the indefinite Toda lattices introduced in
\cite{KO-YE1}.  We  let  for any $\epsilon \in {\mathcal E}$,
\begin{equation}
Z_\epsilon =\left\{\   X=x+ \displaystyle{ \sum_{i=1}^l } (b_i
e_{-\alpha
_i } +  e_{\alpha _i})\ :\ \epsilon_ib_i >0, X\in S({\mathbb R}) \right\}.
\nonumber
\end{equation}
This
is a set of real normalized Jacobi elements, that is, the elements of $Z_{\mathbb
R}$.   Thus the union of all the $Z_\epsilon$ is all of $Z_{\mathbb R}$,
$$ Z_{\mathbb R}=\bigcup_{\epsilon\in{\mathcal E}}Z_{\epsilon}.$$
 The elements
in  $Z_{\mathbb R}$ are thus the signed normalized Jacobi elements (in $S({\mathbb
R})$ and
$Z$ simply denotes $Z_{\epsilon_o}$ where $\epsilon_o=(1,..,1))$.

\begin{defn}
\label{D8.1.2}
{\it Chevalley invariants and isospectral manifold}:
If $x\in \mathfrak g$ and $g_x \in {\mathfrak
g}^\prime$
is defined by $\langle g_x,y \rangle =(x,y)$ for any $y \in {\mathfrak g}$
then the map ${\mathfrak g} \to
{\mathfrak g}^\prime$
sending
$x $ to $g_x$ defines an isomorphism. We can then regard $S({\mathfrak g})$ as the
algebra of polynomial functions on $\mathfrak g$ itself by setting for $f\in
S({\mathfrak
g})$ and $x\in {\mathfrak g}$, $f(x)=f(g_x)$.

The functions $I_1,\cdots,I_l$ now on $\mathfrak g$ and then restricted to $J_{\mathbb R}$,
$Z_{\mathbb R}$ or to $Z_{\mathbb C}$ are called the {\it Chevalley invariants} which are the polynomial functions of $\{a_1,\cdots,a_l,b_1,\cdots,b_l\}$
for $X=\sum_{i=1}^l(a_ih_{\alpha_i}+b_ie_{-\alpha_i}+e_{\alpha_i})$.
The map ${\mathbb I}=(I_1,..,I_l)$ then defines by restriction  a map 
\begin{equation}
{\mathcal I}={\mathcal I}_F: Z_F \to F^l ~.
\nonumber
\end{equation}
Fix $\gamma \in F^l$ in the
image of the map ${\mathbb I}$, and denote 
$$Z(\gamma)_F={\mathcal I}_F^{-1}(\gamma)={\mathbb I}^{-1}(\gamma)\bigcap Z_F,$$
which defines the isospectral manifold of Jacobi elements of $\mathfrak g$.
Note that in the real (isospectral) manifold
$Z(\gamma)$ studied in  \cite{KOS} will just be one
connected component of $Z(\gamma)_{\mathbb R}$.
\end{defn}

\begin{defn}
\label{D8.1.3}
{\it The isotropy subgroup ${\tilde G}^y$ on $\tilde G$}:
 Let $G_{\mathbb C}^y$ be the isotropy subgroup of $G_{\mathbb C}$
 for an element $y\in \mathfrak g_{\mathbb C}$. The group
$G_{\mathbb C}^y$ is an abelian connected algebraic group of complex dimension $l$
(Proposition 2.4 of \cite{KOS}). 
If $y\in {\mathfrak g}$ we
denote by ${\tilde G}^y$ the intersection
\begin{equation}
 {\tilde G}^y=G^y_{\mathbb C} \cap \tilde G .
\nonumber
\end{equation}

If $x\in {\mathfrak g}$ , the centralizer  of $x$ is denoted ${\mathfrak g}^x$ and $dim
{\mathfrak g}^x\geq l$. We say that $x$ is $regular$ if $dim {\mathfrak g}^x=l$.

We consider an open subset of $G_{\mathbb C}$ given by
the biggest piece in the Bruhat decomposition:
\begin{equation}
(G_{\mathbb C})_*=
{\bar N}_{\mathbb C} H_{\mathbb C } N_{\mathbb C} ={\bar N}_{\mathbb C} B_{\mathbb C}~,
\label{8.1.5}
\end{equation}
where $N_{\mathbb C}=\exp ({\mathfrak n}_{\mathbb C}), {\bar N}_{\mathbb
C}=\exp({\bar {\mathfrak n}}_{\mathbb C})$ and $B_{\mathbb C}=H_{\mathbb C}N_{\mathbb C}$.
We let ${\tilde G}_*={\tilde G} \cap (G_{\mathbb C})_*$ and, as in  Notation \ref{D2.2.1}, $B=H_{\mathbb R}\exp ({\mathfrak n})$.
We then have a map
\begin{equation}
 {\bar N}_{\mathbb C}\times  B_{\mathbb C} \to (G_{\mathbb
C})_* ~,
\nonumber
\end{equation}
given by $(n, b)\mapsto nb$ which is an isomorphism of
algebraic varieties. Given $d\in (G_{\mathbb C})_*$, $d$ has a unique
decomposition as
$d=n_{d}b_d$ as in (2.4.6) of \cite{KOS}.
\end{defn}

We caution the reader that ${\tilde G}^y_*$ is not an intersection with $\exp
({\bar
{\mathfrak n}})H\exp({\mathfrak n})$ but rather with  $\exp ({\bar {\mathfrak n}})H_{\mathbb
R}\exp({\mathfrak n})$.
 This is the object that appears, for example, in (3.4.10) of \cite{KOS} and
properly contains (3.2.9) in Lemma 3.2 of  \cite{KOS}.

The following Proposition gives the relation between ${\tilde G}^y$ and $H_{\mathbb R}$:

\begin{prop}
Let $y\in Z_{\epsilon_o}$. Then ${\tilde G}^y$ is
${\tilde G}$ conjugate to the Cartan subgroup $H_{\mathbb R}$ of ${\tilde G}$.\label{P8.1.1}
\end{prop}
\begin{Proof}
By  Lemma 2.1.1 in \cite{KOS} , $y$ is regular. Then as in Lemma 3.2 of \cite{KOS}, $y$
must be conjugate, under an element in $H$,  to an element $x$ in $\mathfrak p$.
Using proposition 2.4 of \cite{KOS}   $G_{\mathbb C}^x$ is connected. Since $x$ is
also  a
regular element it follows that   ${\mathfrak g}_{\mathbb C}^x$ is a Cartan subalgebra and
$G_{\mathbb C}^x$ is a Cartan subgroup of $G_{\mathbb C}$. Using conjugation by an
element
in $K$ we may conjugate this Cartan subgroup if necesary to $H_{\mathbb C}$
(Proposition 6.61 or Lemma 6.62 of \cite{KN}).  We can thus assume that $G^x_{\mathbb
C}=H_{\mathbb C}$. We  obtain that ${\tilde G}^y = {\tilde G} \cap G_{\mathbb C}^y$ is
${\tilde G}$ conjugate to ${\tilde G} \cap H_{\mathbb C}=H_{\mathbb R}$ (see Notation \ref{D2.2.1}).
\end{Proof}

\subsection{Kostant's map $\beta_{\mathbb C}^y$}
 \label{ss8.2} 

Fix $y\in J(\gamma)_{\mathbb R}$. Kostant defines a map 
\begin{equation}
\begin{array}{lll}
 \beta_{\mathbb C}^y ~: ~ &(G_{\mathbb
C}^y)_* & \longrightarrow J(\gamma)_{\mathbb C} \\
& d &\mapsto Ad(n_d^{-1})(y)
\end{array} 
\label{8.2.2}
\end{equation}
with $d=n_db_d, n_d\in \bar N$ and $b_d\in B$.
Note that we have deviated from the convention in \cite{KOS} by exchanging the
roles
of ${\bar N}_{\mathbb C}$ and $N_{\mathbb C}$. We did not exchange the roles of these two
groups in (\ref{8.1.5}) but this is compensated by our use of an inverse in the
definition of the Kostant map.
Theorem 2.4 of \cite{KOS} then implies that $\beta_{\mathbb C}^y$ is an isomorphism of
algebraic varieties.

Denote $\beta^y$ the restriction of $\beta_{\mathbb C}^y$ to the intersection with
$\tilde G$. Thus we have $\beta^y: {\tilde G}^y_* \to Z(\gamma)_{\mathbb R}$ and ${\tilde
G}^y_*=(G^y_{\mathbb C})_*\cap {\tilde G}$

\begin{prop}
Let $y\in J(\gamma)_{\mathbb R}$. The map $\beta^y$
is an isomorphism of smooth manifolds ${\tilde G}_*^y\to J(\gamma)_{\mathbb R}$. \label{P8.2.1}
\end{prop} 

\begin{Proof} The map $\beta^y$ is the restriction to the Lie group
${\tilde G}^y$ of the diffeomorphism of complex analytic manifolds $\beta^y_{\mathbb
C}$.
We obtain that $\beta^y$ must be an injective map. We show surjectivity. If
$z\in
J(\gamma)_{\mathbb R}$ then by surjectivity of $\beta_{\mathbb C}^y$ there is  $g_{\mathbb C}
\in (G_{\mathbb C}^y)_*$ such that $\beta_{\mathbb C}^y (g_{\mathbb C})=z$ and $g_{\mathbb
C}=n_{\mathbb C} b_{\mathbb C}$. Thus $g_{\mathbb C}^c= n_{\mathbb C}^c
b_{\mathbb C}^c$ with $n^c_{\mathbb C}\in {\bar N}_{\mathbb C}$ and $b_{\mathbb C}^c\in B_{\mathbb C}$.
Therefore $\beta^y(g^c_{\mathbb C})= Ad (n_{\mathbb C}^c )^{-1} y$. Since $y^c=y$ we
obtain that $(Ad(n_{\mathbb C})^{-1}y)^c=z^c$. But our assumption is that $z^c=z$.
Hence we have obtained that $\beta_{\mathbb C}^y(g_{\mathbb C}^c)= z$. By the
$injectivity$
of $\beta^y_{\mathbb C}$ we obtain that $g_{\mathbb C}^c= g_{\mathbb C}$. Therefore
$g_{\mathbb C}\in
\tilde G$ and thus $g_{\mathbb C}=g \in { \tilde G}^y$.
This proves $\beta^y$ is a bijection.

By Proposition 2.3.1 of \cite{KOS} , $J(\gamma)_{\mathbb R}$ is a submanifold of real
dimension $l$ of $J(\gamma)_{\mathbb C}$. The diffeomorphism $\beta^y_{\mathbb C}$
restricts to the smooth non-singular map $\beta^y$. Since we have shown that
$\beta^y$ is a bijection, then it is a diffeomorphism.
\end{Proof}

\begin{rem}
\label{R8.2.1}
 If $y\in Z_{\epsilon_o}$ then ${\tilde G}^y$ is a
Cartan subgroup conjugate to $H_{\mathbb R}$ (see Proposition \ref{P8.1.1}). However, in
general it may happen that $y\in J_{\mathbb R}(\gamma)$ is not semisimple or
that  $y$
is semisimple but  ${\tilde G}^y$ is a Cartan subgroup not conjugate to
$H_{\mathbb R}$. For example,
take $y=a_1h_{\alpha_1} + b_1e_{-\alpha_1} + e_{\alpha_1}$ in
the case  of ${\mathfrak{sl}}(2,{\mathbb R})$ ($\tilde G=Ad(SL(2,{\mathbb R})^{\pm})$). We get
a nilpotent
matrix if $a_1^2+b_1=0$. When $a_1^2+b <0$ then $y$ is semisimple but  $\tilde
G^y$ is a compact Cartan subgroup and thus it is not conjugate to $H_{\mathbb
R}$. In
the cases $a_1^2+b >0$ one obtains that ${\tilde G}^y$ is conjugate to
$H_{\mathbb R}$.

 Assume that we are in the case when ${\tilde G}^y$ is conjugate to $H_{\mathbb R}$
(for example $y \in Z_{\epsilon_o}$). Combining Proposition \ref{P8.1.1} and
Proposition \ref{P8.2.1}, Kostant's map gives an imbedding of $Z_{\mathbb R}(\gamma)$,
where $\gamma \in \mathbb R^l$ is in the image of ${\mathbb I}$, into ${\hat H}_{\mathbb
R}$ as an open dense subset.
\end{rem}

\subsection {Kostant's map and toric varieties}
 \label{ss8.3}
We first remark that for
a fixed $y\in Z_{\epsilon_o}(\gamma)$ Kostant's map $\beta^y$
in (\ref{8.2.2}) is just the map $d\to n_d^{-1}$ with $d=n_db_d\in
{\tilde G}^y_*$ and $n_d\in {\bar N}$, so that it
 can be described as a map into the flag manifold: $d\to gB$ in ${\tilde G}/B$,
restricted to ${\tilde G}_*^y$. 

By Proposition \ref{P8.2.1} the map into the flag manifold is a diffeomorphism  onto its image when restricted to ${\tilde G}_*^y$. Hence, since the map is given by the action of a Cartan subgroup on the flag manifold; this action has a trivial isotropy group and the map to the flag manifold sends  ${\tilde G}^y$ diffeomorphically to its image.

  The Cartan subgroup ${\tilde G}^y$  is as good as its conjugate $H_{\mathbb R}$ but, for convenience, we prefer to deal with  $H_{\mathbb R}$ for which we have established notation. We let  $x$ be an element that conjugates
$H_{\mathbb R}$ to ${\tilde G}^y$,  $x^{-1}H_{\mathbb R} x ={\tilde G}^y$. Then the
$\tilde G^y$
orbit of $B$ in ${\tilde G}/B$ is $x^{-1}H_{\mathbb R} x B$. Since we can
translate this set
using multiplication by the fixed element $x$ , we can just study the
$H_{\mathbb R}$
orbit of $xB$ in ${\tilde G}/B$.


We denote:
\begin{itemize}
\item  The map
$q:{\tilde G}^y\to {\tilde G}/B$.
\smallskip

\item  ${\hat Z}(\gamma)= \overline {(q\circ
(\beta^y)^{-1}Z(\gamma)) }=x^{-1}( \overline {(H_{\mathbb R} xB)})$.
\end{itemize}
To study ${\hat Z}(\gamma)$, it is enough to describe in detail the toric variety
$\overline {(H_{\mathbb R}xB)}$. Thus we focus our attention on objects that have this general form:

\begin{defn}
 For any $n\in {\tilde G}$, such that $nB\cap H_{\mathbb R}=\{e\}$, the
toric variety $\overline { (H_{\mathbb R} n B) }$ is called {\it generic} in the sense of \cite{FL-HA1},  if $n\in \underset {w\in W} {\bigcap}w(\bar NB)w^{-1}.$
\end{defn}

We  then assume   $n\in \underset {w\in W} {\bigcap}w(\bar NB)w^{-1}$ and $nB\cap H_{\mathbb R}=\{e\}$ . With these hypotheses we have:
\begin{equation}
nB=\underset {\phi_j\in -\Delta_+}
{\prod}\exp(t^e_je_{\phi_j})B
\nonumber
\end{equation}
 for $t^e_j\in {\mathbb R}$, and 
\begin{equation}
nB=n(w(\Delta))wB=\underset {\phi_j\in -\Delta_+}
{\prod}\exp(t^w_je_{w(\phi_j)})wB
\nonumber
\end{equation}
for  $w \in W$ and $t^w_j\in {\mathbb R}$

We now define:

\begin{defn}
\label{D8.3.2} 
Assume $n\in \underset {w\in W} {\bigcap}w(\bar NB)w^{-1}$ and $nB\cap H_{\mathbb R}=\{e\}$. Denote for any $A\subset \Pi$,
\begin{equation}
n(w(\Delta^A))wB=
\underset{\phi_j\in
-\Delta^A_+}{\prod}\exp(t^w_j e_{w(\phi_j)})wB~.
\nonumber
\end{equation}
We define a map,
\begin{equation}
\begin{array}{llll}
{\check B}~: ~& \quad \quad{\hat H}_{\mathbb R} ~~&\longrightarrow ~~&\overline {( H_{\mathbb R} nB)} \\
&  ( Ad^{w(\Delta^A)}g, [w]^A)&\mapsto & gn(w(\Delta^A))wB
\end{array}
\nonumber
\end{equation}
where $g\in H_{\mathbb R}$ (see Definition \ref{D5.1.1} for ${\hat H}_{\mathbb R}$).
Note here that
\begin{equation}
gn(w(\Delta^A))wB=\underset {\phi_j\in -\Delta_+^A}{\prod}\exp\left(t^w_j
\chi_{w(\phi_j)}(g)e_{w(\phi_j)}\right)wB.
\label{8.3.2}
\end{equation}
\end{defn}
The map $\check B$ can be interpreted as a version of Kostant's map $\beta^y$
for the subsystem determined by the set $A\in \Pi$ and its $w$-translation. A
detailed correspondence with the set ${\overset {\circ}{Z}}_{\mathbb R}
(\gamma)$ could
be made but it requires additional notation.

 Note also  that we have
$\chi_{w(\phi_j)}(g)=\chi_{w(\phi_j)}^{w(\Delta^A)}\circ
Ad^{w(\Delta^A)}(g)$ for $\phi_j\in \Delta^A$. Then we obtain:

\begin{thm}
Assume $n\in \underset{w\in W} {\bigcap}w(\bar NB)w^{-1}$ and $nB\cap H_{\mathbb R}=\{e\}$. The function $\check B$ is a  homeomorphism of topological spaces. The
toric variety $\overline{ (H_{\mathbb R} n B)} $ is a smooth manifold and the map $\check
B$ is a diffeomorphism. \label{T8.3.1}
\end{thm}

\begin{Proof}
First we point out that we already have a smooth manifold $\hat {H}_{\mathbb R}$ 
and the assumption $n\in \underset{w\in W} {\bigcap}w(\bar NB)w^{-1}$ will simply ensure that we can define the map to the flag manifold.

We first show the continuity of $\check B$: We use the
local coordinates, $\{\phi_w:w\in W\}$ for ${\hat H}_{\mathbb R}$ given
in Definition \ref{D4.2.1}. In these local
coordinates, we assume that for each $i=1,..,l$,
 $\phi_{ww_o,i}
(Ad^{ww_o(\Delta^{w_o (A)})}(g), [ww_o]^{w_o(A)})$ (which equals
$\chi^{w(\Delta^A)}_{w(-\alpha_i)}(g)$ or zero) converges to a scalar
$\chi^o_{w(-\alpha_
i )}$.
We note that if $\phi_j=-\displaystyle{\sum_{i=1}^l } c_{i j}\alpha_i $
with each
$c_{ij}$ a non-negative integer,  then
\begin{equation}
 \chi_{w(\phi_j)}=\displaystyle{\prod_{i=1}^l} \chi^{c_{ij}}_{-w(\alpha_i)} =\displaystyle{\prod_{i=1}^l}\phi_{ww_o,i}^{c_{ij}}(Ad^{ww_o(\Delta^{w_o(A)})},
[ww_o]^{w_o(A)}   )~.
\nonumber
\end{equation}
We thus let
\begin{equation}\chi^o_{w(\phi_j)}=\displaystyle{\prod_{i=1}^l}
\left(\chi^{o}_{-w(\alpha_i)}\right)^{c_{ij}}~.
\nonumber
\end{equation}

Let 
\begin{equation}
A^\prime=\{\ \alpha_i\in \Pi: \chi^o_{-w(\alpha_i)}=0 \}.
\nonumber
\end{equation}
Then note that  $\chi^o_{w(\phi_j)}=0$ if and only if $c_{ij}\not=0$ for some
$\alpha_i\in A^\prime$. Thus the only $\chi^o_{w(\phi_j)}$ which are non-zero
correspond to roots in $\Delta^{A^\prime}_+$.

By the assumption made, we have
$$gn(w(\Delta^A))wB=\underset {\phi_j\in -\Delta_+} {\prod}\exp\left(t^w_j
\chi_{w(\phi_j)}(g)e_{w(\phi_j)}\right)wB~,$$
which can be written in terms of the coordinate functions $\phi_{ww_o}$ as
$$\underset {\phi_j\in -\Delta_+}
{\prod} \exp\left(t^w_j  \displaystyle{\prod_{i=1}^l}\phi_{ww_o,i}^{c_{ij}}(Ad^{ww_o(\Delta^{w_o(A)})},[ww_o]^{w_o(A)})
e_{w(\phi_j)} \right) wB~.$$
This then converges (by continuity of $\phi_{ww_o}$) to
$$\underset {\phi_j\in -\Delta_+} {\prod}\exp\left(t^w_j
\chi^o_{w(\phi_j)}(g)e_{w(\phi_j)}\right) wB~,$$
which only involves roots in $\Delta^{A^\prime}_+$, and can be
written as $gn(w(\Delta^{A^\prime}))wB$.
Since any $(Ad^{w(\Delta^A)}(g),[w]^A)\in \hat H_{\mathbb R}$ is completely
determined
by the coordinates $\phi_{w}(Ad^{w(\Delta^A)}(g),[w]^A)$ and some
$(Ad^{w(\Delta^{A^\prime})}(g),[w]^{A^\prime})$ uniquely corresponds to the
coordinates $(\chi^o_{w(\alpha_1)}(g),.., \chi^o_{w(\alpha_l)}(g))$, then
we can
conclude that ${\check B}(Ad^{w(\Delta^A)}(g),[w]^A)$ converges to  $\check B
(Ad^{w(\Delta^{A^\prime})}(g),[w]^{A^\prime})$ whenever the argument
$(Ad^{w(\Delta^A)}(g),[w]^A))$ approaches to $
(Ad^{w(\Delta^{A^\prime})}(g),[w]^{A^\prime})$
in ${\hat H}_{\mathbb R}$. This proves the continuity of the map $\check B$.

Since ${\hat H}_{\mathbb R}$ is compact  (Proposition \ref{P7.3.1}) and ${\check B} ({\hat H}_{\mathbb
R})$ contains the orbit $H_{\mathbb R} n B$ , then  $\overline{ (H_{\mathbb R} n B) }\subset
\check B ({\hat H}_{\mathbb R})$.  From the construction of the map $\check B$ it is
easy to see that its image is contained in  $\overline { (H_{\mathbb R} n B)} $ and thus
${\check B} ({\hat H}_{\mathbb R}) = \overline { (H_{\mathbb R} n B)}$.

The smoothness of the map and that it gives a diffeomorphism follows from the
fact that  $(s_1,..,s_{|\Delta_+|} )\to \displaystyle{\prod_{j=1}^{|\Delta_+| }} \exp ( s_je_{w(\phi_j)} )wB$ constitutes a coordinate system in
the flag
manifold.
\end{Proof}

\begin{rem} We caution the reader that if one replaces the Cartan subgroup   
$H^1_{\mathbb R}$ of $G$ instead of $H_{\mathbb R}$ in the statement of Theorem \ref{T8.3.1} 
then the closure of the orbit $H^1_{\mathbb R}nB$ may not be smooth. In fact the structure of 
 $\overline{ H^1_{\mathbb R}nB }$ can be explicitly described too. Consider only the connected components of each $H_{\mathbb R}(w(\Delta^A)$ associated to $\epsilon$ such that 
$Ad(h_\epsilon) \in Ad( H^1_{\mathbb R})$ in Definition \ref{D4.2.1}. This gives a
subspace of ${\hat H}_{\mathbb R}$ that will correspond to $\overline{ H^1_{\mathbb R}n B}$.
 In terms of the coordinate charts $\phi_w$ one
gets locally ${\mathbb R}^l$ but  now some of  its $2^l$ quadrants may be missing (since  $Ad(H^1_{\mathbb R})$ may have fewer than $2^l$ connected components).  Thus smoothness is obtained exactly when  $Ad(H^1_{\mathbb R})$ contains $2^l$ connected components. Examples that lead to non-smooth closures if one uses the Cartan subgroup $H^1_{\mathbb R}$  are all the $G=SL(n,{\mathbb R} )$ with 
$n$ even.  In terms of the Toda lattice this corresponds to considering the indefinite Toda lattice in \ref{1.1.1} in the Introduction but leaving out some of the signs $\epsilon_i$. When $n=2$, for example,  one obtains a closed interval inside $G/B$ ( which is 
a circle). The disconnected Lie group $\tilde G$ which leads to the Cartan subgroup $H_{\mathbb R}$ is then a requirement in all our constructions and main results.
We thus have:
\begin{cor} Assume $n\in \underset{w\in W} {\bigcap}w(\bar NB)w^{-1}$ and $nB\cap H^1_{\mathbb R}=\{e\}$. Then  $\overline{ H^1_{\mathbb R}n B}$ is smooth
if and only if $Ad(H^1_{\mathbb R})$ has $2^l$ connected components.\label{nonsmooth}
\end{cor}
\label{not-smooth}
\end{rem}

We also have:
\begin{cor} If $n\in \underset{w\in W} {\bigcap}w(\bar NB)w^{-1}$ and $nB\cap H_{\mathbb R}=\{e\}$,   then the toric variety $X=\overline {(H_{\mathbb R} n B)}$ satisfies: $X^{H_{\mathbb R}}={(G/B)}^{H_{\mathbb R} }$, the $ H_{\mathbb R}$ fixed points. \label{genericity2}
\end{cor}
\begin{proof} We use Theorem \ref{T8.3.1}. The manifold $\hat{H}_{\mathbb R}$ has an $H_{\mathbb R}$ action and the only fixed points are the $(e, w)$ with $w\in W$. These get mapped to the fixed points  of the $H_{\mathbb R}$ action in $G/B$, the cosets $wB$, $w\in W$.  Thus $X^{H_{\mathbb R}}={(G/B)}^{H_{\mathbb R}}$.  This also follows  directly if 
if we consider  \ref{8.3.2} with $A=\emptyset$ and we just let each $\chi_{w(\phi_j)} $ go to zero and obtain $wB$.
\end{proof}

\begin{rem}
\label{R8.3.1} $\text{}$
\begin{enumerate}
 \item By Theorem 3.6 in  \cite{FL-HA2}  and Remark 3 in
p. 257 of \cite{FL-HA1} we may assume that  the map into the flag manifold 
which appears as a consequence of Kostant's map in Subsection \ref{ss8.3}  is such that one in fact obtains the $H_{\mathbb R}$ orbit of an element  $n\in \underset
{w\in W} {\bigcap}w(\bar NB)w^{-1}$.

\item In \S 7 of  \cite{CA-KU} and in \cite{DA}  a more restrictive definition of the notion of {\em genericity }  is given. This discussion is only relevant  for the case of toric varieties in $G/P$ where
$P$ is a parabolic subgroup rather than just a Borel subgroup as is the case in the present work. In these more general situations sometimes generic varieties as defined in \cite{FL-HA1} are not normal. In any case, note that  if $P=B$ our Corollary \ref{genericity2} implies that if $n\in \underset{w\in W} {\bigcap}w(\bar NB)w^{-1}$  and $nB\cap H_{\mathbb R}=\{e\}$ then the corresponding toric variety is also generic in the sense discussed in  \cite{CA-KU} or \cite{DA}.  Finally we observe that normality is not used  in any of our results.  Here we rely instead on explicit coordinate charts to obtain the smoothness of our toric varieties and describe their topological structure. (We would like to thank H. Flaschka for sending the papers \cite{CA-KU,DA} to us.)
\end{enumerate}
\end{rem}
We conclude:
\begin{thm}
Let $\gamma\in {\mathbb R}^l$ , then ${\hat
Z}(\gamma)_{\mathbb R}$ is a smooth compact manifold diffeomorphic to ${\hat
H}_{\mathbb R}$.  \label{T8.3.2}
\end{thm}
\begin{Proof} This is just Theorem \ref{T8.3.1} and  the definition of ${\hat
Z}(\gamma)$. The two conditions in Theorem  \ref{T8.3.1} are satisfied by  
Proposition \ref{P8.2.1} (Subsection \ref{ss8.3})  and  by Theorem 3.6 in  
\cite{FL-HA2} and Remark 3 in p. 257 of \cite{FL-HA1} as noted above in Remark \ref{R8.3.1}.
\end{Proof}

\bibliographystyle{amsplain}

\end{document}